  \documentclass[11pt]{article}       
  \usepackage[cp1251]{inputenc}               
  \usepackage[english,russian]{babel}         
  \usepackage{literat}                        
  \usepackage{amsfonts,amssymb,amsmath,
              amstext,amsbsy,amsopn,amscd}    
  \usepackage{amsthm}                         
  \usepackage[pagestyles]{titlesec}           
  \usepackage[labelsep=period,
              font=footnotesize]{caption}     
  \usepackage{graphicx}                       

  \usepackage[all,knot]{xy} \xyoption{arc} 
  \usepackage{eepic} 
  \def\ChangeFontdimen#1#2#3#4#5{{#1%
     \setlength{\fontdimen2\font}{#2}%
     \setlength{\fontdimen3\font}{#3}%
     \setlength{\fontdimen4\font}{#4}%
     \setlength{\fontdimen7\font}{#5}}}
  \mag=\magstep2
\advance\marginparwidth-28.mm \reversemarginpar 
\newtitlemark{\AUTHORmark} \def\AUTHORmark{\relax}
\newtitlemark{\TITLEmark}  \def\TITLEmark{\relax}
\newpagestyle{sbornik}{\headrule
               \sethead[\fontsize{9}{0}\selectfont\thepage][\AUTHORmark][]%
                       {}{\TITLEmark}{\fontsize{9}{0}\selectfont\thepage}}
\def\TITLE#1#2#3#4#5{%
  \clearpage\thispagestyle{empty}%
  \mark{}%
  \def\AUTHORmark{\fontsize{7}{0}\selectfont\spaceskip 0.19em #1}%
  \def\TITLEmark{\fontsize{9}{0}\textsc{#3}}%
  \def\thefootnote{\fnsymbol{footnote}}%
  \begin{center}%
     \baselineskip14.0pt
     \vspace*{20.0mm}%
     {\fontsize{9}{0}\selectfont\spaceskip 0.17em #1}%
     \par\vspace{5.0mm}%
     {\fontseries{b}\selectfont\spaceskip 0.51em #2}%
     \ifcase#4\or\footnotetext{\noindent#5}\fi
     \par\vspace{1.2mm}%
  \end{center}\par}
  \makeatletter\renewcommand{\@biblabel}[1]{#1.}\makeatother
  \def\REFERENCES#1{\begin{thebibliography}{#1}\fontsize{9}{0}\selectfont
     \parskip.8ptplus0.2pt\itemsep0.pt}

  \def\Jr{,\,Jr.}
  \def\ifxtmp#1#2{\def\ttt{#1}\ifx\ttt\tmp\def\next{#2}\fi}
  \def\razryadka#1{\def\tmp{#1}%
        \def\next{ #1\razryadka}%
        \ifxtmp +   {\razryadka}%
        \ifxtmp .   {.\razryadka}%
        \ifxtmp ,   {,\hskip3.3ptplus1.8ptminus.1pt\razryadka}%
        \ifxtmp ~   {\hskip2.9ptplus.3ptminus.1pt\razryadka}%
        \ifxtmp |   {\hskip3.5ptplus1.ptminus.5pt}%
        \ifxtmp \Jr {,\hskip2.5ptplus.2ptminus.1ptJr.\razryadka}%
      \next}
  \titlespacing {\section}
                {0.pt} {3.explus.5exminus.3ex} {2.3explus.6exminus.2ex}
  \titleformat  {\section}            [hang]
        {\fontsize{9}{11}\fontseries{b}\selectfont
             \spaceskip 0.4em\filcenter}  {\Rus{\S}\,\thesection.}{.7em}{}

  \titlespacing {\subsubsection}
                {\parindent} {1.5explus.1ex} {2.5mm}
  \titleformat  {\subsubsection}      [runin]
        {\normalfont\bfseries}
        {\thesubsubsection.}{.6em}{}[.\hspace{1.5ptplus.9pt}\null] 
  \titlespacing {\subparagraph}
                {\parindent} {1.5explus.1ex} {2.5mm}
  \titleformat  {\subparagraph}       [runin]
        {\normalfont\bfseries}
        {\thesubparagraph.}{}{}[.\hspace{1.5ptplus.9pt}\null]   
    \newskip\prethm \prethm 2.7pt plus1.1pt minus0.3pt
    \newskip\posthm \posthm 2.7pt plus1.1pt minus0.3pt
\makeatletter
  \def\thmhead@plain#1#2#3{%
    \thmname{#1}\thmnumber{\@ifnotempty{#1}{ }\@upn{#2}}%
    \thmnote{ \hskip0.4emplus0.15emminus0.1em{\the\thm@notefont#3}}}
\makeatother
    \newdimen \sVrazr \sVrazr 0.28em                    
    \newskip  \pVrazr \pVrazr 0.7em plus.2em minus.1em  
    \newtheoremstyle{STATEMENT}%
       {\prethm}{\posthm}{\itshape}{\parindent}{\spaceskip\sVrazr}%
         {.}{\pVrazr}{}
    \newtheoremstyle{STATEMENTx}
       {\prethm}{\posthm}{\itshape}{\parindent}{}%
         {.}{\pVrazr}{\thmnote{#3}}
    \newtheoremstyle{EXPLANATION}%
       {\prethm}{\posthm}{}{\parindent}{\spaceskip\sVrazr}{.}{\pVrazr}{}
    \newtheoremstyle{EXPLANATIONx}
       {\prethm}{\posthm}{}{\parindent}{}%
         {.}{\pVrazr}{\thmnote{#3}}
     \theoremstyle{STATEMENT}
  \newtheorem  {unLem}  {Л е м м а             \,} 
  \newtheorem  {unCor}  {С л е д с т в и е     \,} 
  \newtheorem  {unThe}  {Т е о р е м а         \,} 
     \theoremstyle{STATEMENTx}
  \newtheorem* {itXxx}  {}              
     \theoremstyle{EXPLANATION}
  \newtheorem* {nnExa}  {П р и м е р}              
  \newtheorem  {unRem}  {З а м е ч а н и е     \,} 
  \newtheorem  {unDef}  {О п р е д е л е н и е \,} 
     \theoremstyle{EXPLANATIONx}
  \newtheorem* {rmXxx}  {}              
\newcounter{NNN}
\newcounter{nnn}[NNN]
\renewcommand{\thennn}{\theNNN.\arabic{nnn}}
\def\Nit#1#2{\refstepcounter{nnn}\begin{itXxx}%
             [\thennn\ifcase#1\or
              .\hskip0.7em plus.2em minus.1em{\spaceskip\sVrazr #2}\fi]}
\def\Nrm#1#2{\refstepcounter{nnn}\begin{rmXxx}%
             [\thennn\ifcase#1\or
              .\hskip0.7em plus.2em minus.1em{\spaceskip\sVrazr #2}\fi]}
\def\endNit{\end{itXxx}}
\def\endNrm{\end{rmXxx}}
 \newskip\preprf \preprf\prethm
 \newskip\posprf \posprf\posthm
 \newskip\SAVE
 \def\PROOF  {\vskip\preprf{\spaceskip\sVrazr
                          Д о к а з а т е л ь с т в о}\hskip\pVrazr}
 \def\PROOFp {\vskip\preprf{\spaceskip\sVrazr
                          Д о к а з а т е л ь с т в о}.\hskip\pVrazr}

 \def\ENDPROOF {\vskip\posprf}

   \DeclareMathOperator {\ind}   {ind}

\begin{document}
%
  \sfcode`\.=3000 \sfcode`\?=3000\sfcode`\!=3000
  \frenchspacing
\let\savess\thesubsection
\def\THESUBSECTION{\savess}
%
 \def\paperRus{\Russian
  \ChangeFontdimen{\normalfont}          {0.37em} {0.13em} {0.10em} {0.13em}%
  \ChangeFontdimen{\normalfont\itshape}  {0.38em} {0.12em} {0.09em} {0.13em}%
  \ChangeFontdimen{\normalfont\bfseries} {0.39em} {0.14em} {0.11em} {0.13em}%
       \setcounter{section}  0
       \setcounter{equation} 0
   \renewcommand{\theequation}{\arabic{equation}} 
   \renewcommand{\thesubsection}{\THESUBSECTION}
       \setcounter{unCor}    0
       \setcounter{unDef}    0
       \setcounter{unExa}    0
       \setcounter{unLem}    0
       \setcounter{unPro}    0
       \setcounter{unRem}    0
       \setcounter{unThe}    0
       \setcounter{figure}   0
       \setcounter{NNN}      0
           \renewcommand{\refname}%
     {\centerline{\fontsize{9}{0}\textmd{Л\,И\,Т\,Е\,Р\,А\,Т\,У\,Р\,А}}}%
           \renewcommand{\figurename}{Рис.}%
               \input}
%
     \paperRus  
%


\newcommand{\ZZ}{\mathbb Z}
\newcommand{\RR}{\mathbb R}
\newcommand{\geom}{\mathrm {geom}}
\renewcommand{\int}{{\mathrm{int\,}}}
\newcommand{\sgn}{{\mathrm{sgn\,}}}

\TITLE {S.\, A.\, B O G A T Y I, \, \, E.\, A.\, K U D R Y A V T S E V A,
\, \, H.\, Z I E S C H A N G} {ON INTERSECTIONS OF CLOSED CURVES ON SURFACES} 
{} {1} {This research was done during the visit and working of the third author 
at the Department of Mathematics and Mechanics of Moscow State M.\,V.~Lomonosov University in 2001--2003. This visit was supported by the Stiftungsinitiative Johann Gottfried Herder. The first and second authors were supported by the Russian Foundation for Basic Research (grant nos.\ 09--01--00741-a and 10--01--00748-a), the Programme of the President of the Russian Federation for the Support of Leading Scientific Schools (grant nos.\ 
НШ-1562.2008.1 and НШ-3224.2010.1, the programme ``Development of Scientific Potential in Higher Education'' (grant no.\ 2.1.1.3704), and the federal target programme ``Scientific and Scientific-Pedagogical Personnel of Innovative Russia'' (grant no.\ 14.740.11.0794).}

\def\mapr#1{\smash{\mathop{\buildrel{#1}\over\longrightarrow}}}
\def\mapd#1{\Big\downarrow\rlap{$\vcenter{\hbox{$#1$}}$}}
\def\mapu#1{\Big\uparrow\rlap{$\vcenter{\hbox{$#1$}}$}}


{\small
The problem on the minimal number (with respect to deformation) of intersection points of two closed curves on a surface is solved. Following the Nielsen approach, we define classes of intersection points and essential classes of intersection points, which ``are preserved under deformation'' and whose total number is called the Nielsen number. If each Nielsen class consists of a unique point and has a non-vanishing index after a suitable deformation of the pair of curves, one says that {\it the Wecken property holds}. We compute the minimal number of intersection points in terms of the Nielsen numbers and the Reidemeister numbers. In particular, we prove that the Wecken property does not hold for some pairs of closed curves. Moreover, all the non-vanishing indices of the Nielsen classes equal $\pm1$, while the non-vanishing Jezierski semi-indices equal $1$. Similar questions are studied for the self-intersection problem of a curve on a surface.
}



\begin{center} 
{\small
С.\, А.\, БОГАТЫЙ, \, \, Е.\, А.\, КУДРЯВЦЕВА, \, \, Х.\, ЦИШАНГ} \bigskip \\ \bf 
О ПЕРЕСЕЧЕНИЯХ ЗАМКНУТЫХ КРИВЫХ \\ НА ПОВЕРХНОСТЯХ \bigskip
\end{center}

\begin{abstract}
В работе решена задача о минимальном (относительно деформации) числе
точек пересечения двух замкнутых кривых на поверхности. Следуя
Нильсену, определяются классы точек пересечения и существенные классы
точек пересечения, которые ``сохраняются при деформации'' и число
которых называется числом Нильсена. Если все классы Нильсена состоят
из одной точки и имеют ненулевой индекс после подходящей деформации
пары кривых, говорят, что {\it выполнено свойство Векена}. Вычислено
минимальное число точек пересечения в терминах чисел Нильсена и чисел
Райдемайстера. В частности, мы доказываем, что свойство Векена
выполнено не для любых пар замкнутых кривых. Более того, отличные от
нуля индексы классов Нильсена равны $\pm1$, а отличные от нуля
полуиндексы Езерского равны $1$. Такие же вопросы изучены для задачи
о самопересечении кривой на поверхности.
\end{abstract}

\section{ВВЕДЕНИЕ} \label{sec:1}
Изучение задачи о совпадении для пары отображений двумерного тора в
поверхность~\cite {BKZ} привело авторов к исследованию свойства
Векена в задачах о пересечении и самопересечении замкнутых кривых на
поверхностях, а также к исследованию некоторых топологических свойств
(существенность, ``специальность'', ``тривиальность'') классов
Нильсена точек пересечения и самопересечения. В данной работе мы
полностью решаем задачи о пересечении и самопересечении замкнутых
кривых на поверхностях. Ответы даются в терминах чисел Нильсена и
чисел Райдемайстера задач о пересечении и самопересечении. В
частности, мы показываем, что свойство Векена выполнено не для любых
замкнутых кривых и пар замкнутых кривых.

Опишем {\it задачу о пересечении}. Даны два непрерывных
отобра\-же\-ния $f_1\colon\, X_1 \to Y$, $f_2\colon\, X_2 \to Y$
между топологическими пространствами. Возникает задача о нахождении
представителей гомотопических классов этих отображений, у которых
число точек пересечения минимально среди всех пар представителей
классов $[f_1]$ и $[f_2]$. В частности, нужно найти {\it минимальное
число точек пересечения гомотопических классов} $[f_1]$ и $[f_2]$:
 $$
 MI[f_1,f_2] := \min_{f'_1 \simeq f_1,f'_2 \simeq f_2}|\int(f'_1,f'_2)|,
 $$
где $\simeq$ означает гомотопность отображений,
 $$
 \int(f_1,f_2) := \{(x_1,x_2)\in X_1\times X_2 \mid f_1(x_1) = f_2(x_2)\} .
 $$
Подчеркнем, что под {\it точкой пересечения}
$(x_1,x_2)\in\int(f_1,f_2)$ мы понимаем точку декартова произведения
$X_1\times X_2$, а не ее проекции на сомножители $X_1$ и $X_2$ и не
их образ в $Y$. Для решения этой задачи, следуя подходу
Я.~Нильсена~\cite {Ni}, разумно разбить множество точек пересечения
отображений $f_1$ и $f_2$ на классы: две точки $(x_1,x_2)$,
$(x_1',x_2')$ пересечения принадлежат одному {\it классу Нильсена},
если существует пара путей $\alpha_i\colon\,[0,1]\to X_i$, $i=1,2$,
со свойствами
 $$
\alpha_i(0) = x_i,\quad \alpha_i(1) = x'_i, \qquad i=1,2, \qquad \quad f_1
\circ \alpha_1 \simeq_\partial  f_2 \circ \alpha_2 .
 $$
Здесь $\simeq_\partial$ означает, что при гомотопии образы концов
остаются на месте.

Если $X_1,X_2,Y$ являются многообразиями размерностей $\dim X_i=n_i$
и $\dim Y=n=n_1+n_2$, можно вести {\it индекс} изолированной точки
пересечения. Например, для изолированной {\it точки пересечения}
$(x_1^0,x_2^0)\in\int(f_1,f_2)$ мы возьмем окрестности
$U_{X_1},U_{X_2}$ и $U_Y$ точек $x_1^0$, $x_2^0$ и
$y^0=f_1(x_1^0)=f_2(x_2^0)$ соответственно, так, что $(x_1^0,x_2^0)$
является единственной точкой пересечения в $U_{X_1}\times U_{X_2}$ и
$f_i(U_{X_i})\subset U_Y$. Эти три окрестности отождествляются с
единичными шарами $D^{n_i}\subset\RR^{n_i}$ и $D^n\subset\RR^n$ при
помощи (сохраняющих ориентацию, если $X_1$, $X_2$ и $Y$
ориентированы) гомеоморфизмов $U_{X_i}\to D^{n_i}$ и $U_Y\to D^n$,
переводящих точки $x_i^0$ и $y^0$ в центры шаров, $i=1,2$. Степень
отображения
 \begin{equation}\label{Gamma}
 \Gamma\colon\, \partial(D^{n_1}\times D^{n_2}) \cong S^{n-1}\to S^{n-1},\quad
 (x_1,x_2) \mapsto \frac{f_1(x_1) - f_2(x_2)}{|f_1(x_1) - f_2(x_2)|}
 \end{equation}
называется {\it индексом $\ind_{(x_1^0,x_2^0)}(f_1,f_2)$ точки
пересечения} $(x_1^0,x_2^0)$; она не зависит от выбора окрестностей и
меняет знак, если меняется ориентация в $U_{X_1}$ или в $U_{X_2}$ или
в $U_Y$. (Если $X_1$ или $X_2$ или $Y$ неориентируемо, индекс
определяется по модулю $2$.) Индекс можно определить для множеств
точек пересечения; например, индекс конечного множества изолированных
точек пересечения равен сумме индексов точек этого множества.
Ю.~Езерский~\cite {DobJez,Jez} определил {\it полуиндекс класса
Нильсена} (для задачи о совпадении), см.\ определение~\ref
{def:2.1}~(d,~e), который принимает значения в $\ZZ_+$ и является
более тонким инвариантом, чем индекс, если $X_1$ или $X_2$ или $Y$
неориентируемо. Если $X_1,X_2$ и $Y$ ориентируемы, полуиндекс класса
равен абсолютному значению индекса этого класса; если $X_1$ или $X_2$
или $Y$ неориентируемо, индекс получается из полуиндекса
проектированием $\ZZ_+$ в $\ZZ_2$.

Класс Нильсена точек пересечения называется {\it существенным}, если
его индекс не нуль. Число $NI[f_1,f_2]$ существенных классов точек
пересечения называется {\it числом Нильсена} пары $f_1$ и $f_2$ для
задачи о пересечении. Существенный класс точек пересечения
``сохраняется при гомотопии'', не меняя индекс. Поэтому число
Нильсена $NI[f_1,f_2]$ является инвариантом пары гомотопических
классов $[f_1]$, $[f_2]$. Ясно, что $MI[f_1,f_2]\ge NI[f_1,f_1]$.
Если числа $MI[f_1,f_2]$ и $NI[f_1,f_2]$ совпадают, говорят, что
выполняется {\it свойство Векена}~\cite {We} для задачи о пересечении
пары $f_1$ и $f_2$. В этом случае существуют отображения $f'_1\simeq
f_1$ и $f'_2\simeq f_2$, для которых каждый класс Нильсена точек
пересечения является существенным и содержит ровно одну точку.
Отметим, что рассматриваемая во многих работах~\cite{BGKZ, BGKZ1,
GKZ} задача о корнях для отображения $f_1\colon\, X_1 \to Y$ между
многообразиями совпадает с задачей пересечения для пары отображений
$(f_1,f_2)$, где $f_2\colon\, X_2\to Y$ и $X_2$ одноточечно.

Важная теорема Добренько и Кухарского~\cite{Kukh} (см.\ также~\cite{McCord}) состоит в том,
что для ориентируемых компактных многообразий размерностей $\dim X_1=n_1$, $\dim
X_2=n_2$ и $\dim Y=n=n_1+n_2$ любая пара непрерывных отображений
обладает свойством Векена для задачи о пересечении в предположении
$\max\{n_1,n_2\}\ne2$. Последнее предположение нельзя отбросить, как
показало решение задачи о корнях в размерности два~\cite {BGKZ,
BGKZ1, GKZ}, равносильной задаче о пересечении в размерностях
$(n_1,n_2)=(2,0)$. В размерностях $n_1=n_2=2$ свойство Векена
выполнено и доказывается аналогично случаю $n_1,n_2\ge3$ при помощи
аналога трюка Уитни \cite {Kervaire,Wall,Milnor,Kukh,McCord}, при
этом каждый двумерный ``диск Уитни'' имеет, вообще говоря, точки
самопересечения и точки пересечения с $f_1(X_1)$ и $f_2(X_2)$, однако
от этих пересечений можно избавиться при помощи деформации
отображений $f_1$ и $f_2$ посредством ``пальцевых движений'' Кассона
(при которой не возникают новые точки пересечения $f_1$ и $f_2$, а
лишь точки самопересечения). Интересно отметить, что для {\it задачи
о совпадении} отображений между поверхностями, тесно связанной с
задачей о пересечении в размерностях $n_1=n_2=2$, свойство Векена не
всегда выполнено~\cite {Bo1,Bo2,Ho,BGKZ1}, и здесь выявление пар со
свойством Векена удалось только в специальных случаях, например,
когда отображения являются гомеоморфизмами~\cite {Bo0,Iv}, или когда
$Y$ является сферой, проективной плоскостью или тором~\cite
{Jez1,Jez2,BBPT}, или когда одно из отображений является константным
(задача о корнях)~\cite{BGKZ, BGKZ1, GKZ}, или когда $X$ является
тором~\cite{BKZ}.

С задачей о пересечении тесно связана {\it задача о
самопересе\-че\-нии}. Для непрерывного отображения $f\colon X \to Y$
между тополо\-ги\-чес\-кими пространствами {\it множеством
самопересечения} назовем множество $\int(f) := \{(x_1,x_2)\in
(X\times X)\setminus\Delta \mid f(x_1) = f(x_2)\}$, где
$\Delta:=\{(x,x)\mid x\in X\}$ --- диагональ. Аналогично возникает
задача о {\it минимальном числе точек самопересечения гомотопического
класса} $[f]$:
 $$
 MI[f] := \min_{f' \simeq f} |\int(f')|.
 $$
Опять для изучения этой задачи разумно разбить точки
самопе\-ре\-се\-че\-ния отображения $f$ на классы: две точки
$(x_1,x_2)$, $(x_1',x_2')$ самопересечения принадлежат одному {\it
классу Нильсена}, если суще\-ст\-вует пара путей
$\alpha_1,\alpha_2\colon\,[0,1]\to X$ со свойствами
 $$
\alpha_i(0) = x_i,\quad \alpha_i(1) = x'_i, \qquad i=1,2, \qquad \quad f
\circ \alpha_1 \simeq_\partial  f \circ \alpha_2 .
 $$
Понятия индекса точки самопересечения, существенного класса Нильсена,
числа $NI[f]$ и свойства Векена для отображения замк\-ну\-того
$n$-мерного многообразия $X$ в $2n$-мерное многообразие $Y$
определяются как в задаче о пересечении. Подчеркнем, что {\it точка
самопересечения} $(x_1,x_2)\in\int(f)$ является упорядоченной парой
$(x_1,x_2)$, т.е.\ элементом пространства $(X\times
X)\setminus\Delta$. Рассмотрение неупорядоченной пары $\{x_1,x_2\}$,
т.е.\ элемента пространства $(\exp_2X)\setminus X\cong((X\times
X)\setminus\Delta)/(x_1,x_2)\sim(x_2,x_1)$, дает другое понятие точки
самопересечения. В настоящей работе такая точка самопересечения
$\{x_1,x_2\}$ называется {\it геометрической точкой самопересечения},
а минимальное число геометрических точек самопересечения
гомо\-то\-пи\-чес\-кого класса $[f]$ обозначается через
$MI_\geom[f]$. Для гео\-мет\-ри\-чес\-ких точек самопересечения
аналогично определяются экви\-ва\-лент\-ность по Нильсену, индекс,
существенный класс Нильсена, число Нильсена $NI_\geom[f]$ и свойство
Векена.

Подчеркнем, что $\int(f)\ne\int(f,f)$, а лишь
$\int(f)\subset\int(f,f)$. Кроме того, в задаче о пересечении для
пары отображений $(f,f)$ отображения $f_1=f$ и $f_2=f$ деформируются
независимо, поэтому нет очевидной связи между числами $MI[f]$ и
$MI[f,f]$, а также между числами $NI[f]$ и $NI[f,f]$.

Уитни доказал, что любое $n$-мерное замкнутое многообразие вложимо в
$\RR^{2n}$, и предложил метод, позволяющий доказать свойство Векена
$MI[f]=NI[f]$ и $MI_\geom[f]=NI_\geom[f]$ в задаче о самопересечении
для отображения $n$-мерного замкнутого много\-об\-ра\-зия $X$ в
$2n$-мерное многообразие $Y$ при $n\ge3$~\cite {Whitney} (см.\ также
\cite {Wall,Milnor}). Этот метод называется {\it трюком Уитни} и
применяется во многих работах~\cite {Sark}. При $n=2$ существуют
отображения, для которых свойство Векена не выполнено~\cite{KeMi}.

В теореме~\ref{pro:2.4} и следствии~\ref {cor:geomselfint} мы решаем
задачу о самопересечении для замкнутой кривой на поверхности (случай
$n = 1$) и даем критерий справедливости свойства Векена для заданной
кривой. Мы также доказываем, что отличные от нуля полуиндексы классов
Нильсена точек самопересечения (и геометрических точек
самопересечения) равны $1$, и находим минимальное число классов
Нильсена, являю\-щих\-ся существенными и/или специальными и/или
тривиальными. Оказывается, что на связной поверхности выполнено
свойство Векена для любой замкнутой кривой в том и только том случае,
когда эта поверхность гомеоморфна плоскости, сфере или проективной
плоскости (т.е.\ имеет конечную фундаментальную группу). Таким
образом, свойство Векена для задачи о самопересечении при $n=1$ не
всегда справедливо.

В теореме \ref {pro:2.2} мы решаем задачу о пересечении для пары
замкнутых кривых на поверхности (т.е.\ в размерностях $n_1=n_2=1$) и
даем критерий справедливости свойства Векена для такой пары кривых.
Мы также доказываем, что отличные от нуля полуиндексы классов
Нильсена точек пересечения равны $1$, и находим минимальное число
классов Нильсена, являющихся существенными и/или специальными.
Оказывается, что на связной поверх\-ности выполнено свойство Векена для любой
пары замкнутых кривых в том и только том случае, когда эта
поверхность ориентируема или гомеоморфна $\RR P^2$. В частности, на любой неориентируемой связной поверхности, отличной от $\RR P^2$, имеются пары кривых, которые не обладают свойством Векена
для задачи о пересечении.

С помощью теорем~\ref {pro:2.4}(a) и \ref {pro:2.2} настоящей работы,
включая утвер\-жде\-ния о существенных, специальных и тривиальном
классах Ниль\-сена точек самопересечения и пересечения, авторами
решена задача о совпадении для пар отображений двумерного тора в
поверхность~\cite[теоремы~1.1 и~3.6]{BKZ}. Методы, используемые в
настоящей работе, основаны на наблюдении
Пуанкаре~\cite[с.~465-475]{Po}, состоящем в том, что на компактной
связной поверхности неположительной эйлеровой характеристики
минимальное число точек (само)пересечения имеют кривые, являющиеся
кратчайшими относительно гиперболической римановой метрики. Часть
результатов настоящей статьи (см.\ тео\-ре\-мы~\ref {pro:2.4}~(a)
и~\ref {pro:2.2}~(a)) была получена Тураевым и Виро~\cite {TuVi} с
использованием этого наблюдения. Для замкнутых кривых на
ориен\-ти\-ру\-емых поверхностях задачи о самопересечении и
пересечении решены в работах~\cite{Ca,CaCh} в терминах классов
эквивалентности точек (само)пересечения, целиком состоящих либо из
классов Нильсена (для задачи о пересечении), либо из классов строгой
эквивалентности по Нильсену (для задачи о самопересечении), где
индексы классов принимают значения в $\ZZ$. Алгоритмы, выясняющие,
можно ли продеформировать данные замкнутые кривые на поверхности в
кривые без (само)пересечений, имеются в~\cite
{Rein,Z,Ch,BiSe,Lu1,Lu2}. В задачах о пересечении и самопересечении
незамкнутых кривых с закрепленными концами на крае поверхности
свойство Векена было доказано Тураевым~\cite{Turaev}, см.\
также~\cite{GKZ2,GKZ3}.

Результаты настоящей статьи и упомянутые выше результаты
по\-ка\-зы\-вают, что в задаче о пересечении пары отображений
$f_i\colon\,X_i^{n_i}\to Y^{n_1+n_2}$, $i=1,2$, $n_1\le n_2$, между замкнутыми
многообразиями открытым остается только случай размерностей
$(n_1,n_2)=(1,2)$, а в задаче о самопересечении
отображения $f\colon\,X^n\to Y^{2n}$ между многообразиями -- только
``двумерный'' случай $n=2$ (изучаемый, например, в~\cite {ScTe}).
По-видимому, нерешенными являются также задачи о {\it минимальном
числе точек самопересечения регулярных кривых} на поверхности $S$
среди регулярных кривых данного класса регулярной гомотопии, а также
о нижней оценке числа самопересечения для нестягиваемых кривых,
гомотопические классы которых принадлежат фиксированной подгруппе
$H\subset\pi_1(S)$, например, $H=\langle\!\langle g\rangle\!\rangle$,
нормальное замыкание элемента $g\in\pi_1(S)$.

Статья имеет следующее строение. В \S\ref {sec:2} даются основные
определения задач пересечения и самопересечения в ``одномерном''
случае, а также мотивируется введение нами более тонких понятий, чем
приведенные выше. В \S\ref {sec:3} решаются задачи о минимальном
числе точек (само)пе\-ре\-се\-че\-ния для кривых на цилиндре и листе
Мебиуса (леммы~\ref {lem:Cylinder}, \ref {lem:Moebius} и
следствия~\ref {cor:geomCylinder}, \ref {cor:geomMoebius}). В \S\ref
{sec:4} вводятся понятия специальных кривых и специальных пар кривых.
В \S\ref {sec:5} и \ref {sec:6} решаются задачи о минимальном числе
точек самопересечения и пересечения для замкнутых кривых на
поверхностях (теоремы \ref {pro:2.4}, \ref {pro:2.2} и следствие~\ref
{cor:geomselfint}).

Авторы благодарны Д.Л.\ Гонсалвесу за указание на работу~\cite {Jez}
и У.~Кошорке за указание определения ``строго эквивалентных'' по
Нильсену точек самопересечения кривых, см.\ определение~\ref
{def:2.1}~(b).

\section{ОСНОВНЫЕ ОПРЕДЕЛЕНИЯ}\label{sec:2}

Основные понятия, введенные в данном разделе --- это описанные
в~\S\ref {sec:1} понятия {\it эквивалентность по Нильсену} точек
пересечения или самопересечения, {\it индекс} класса Нильсена, {\it
существенный класс Нильсена}, {\it число Нильсена} и {\it свойство
Векена}. Этих понятий достаточно для полного решения задач о
пересечении и самопересечении замкнутых кривых на поверхности. Другие
(более тонкие) понятия и инварианты не являются основными и введены
нами для изучения следующих вопросов.

1) В работе~\cite{BKZ} использованы без доказательства некоторые
топо\-ло\-ги\-чес\-кие свойства классов Нильсена точек пересечения и
само\-пе\-ре\-се\-чения. Для до\-ка\-за\-тель\-ства этих свойств мы
вводим и изучаем сле\-дующие классы Нильсена: (i) {\it специальные
классы Нильсена} точек (само)пере\-се\-че\-ния; (ii) {\it тривиальный
класс} Нильсена точек самопересечения.

2) Как показывают полученные нами ответы в задачах о
само\-пе\-ре\-се\-чении и пересечении (теоремы \ref {pro:2.4} и \ref
{pro:2.2}), свойство Векена невы\-пол\-нено для некоторых замкнутых
кривых и пар замкнутых кривых. Возникает задача об изучении
тех классов Нильсена точек самопересечения или пересечения, которые содержат пары эквивалентных точек с противоположными индексами. Мы изучаем топологические свойства этих классов: показываем, что они являются (строго) специальными, см.\ определение~\ref {def:2.1}~(c), и в случае точек самопересечения исследуем, сколько из этих классов являются (строго) геометрически специальными. Мы также изучаем геометрические свойства этих классов: вычисляем их индексы, полуиндексы и наименьшее число точек в каждом классе; выясняем, какие из них являются существенными, а какие (строго) дефективными.

3) Так как свойство Векена невыполнено (см.\ выше), 
представляется естественной задача о введении более тонкого отношения эквивалентности на
множестве точек (само)пе\-ресечения, чем эквивалентность по Нильсену,
а также более тонкого понятия индекса класса, и об изучении выполнения
аналога свойства Векена для полученных более тонких инвариантов.
В связи с этим мы вводим следующие понятия: (i) {\it строго
эквива\-лентные по Нильсену} точки самопересечения; 
(ii) {\it полуиндекс} класса Нильсена точки (само)пересечения 
и необходимые для него понятия {\it самосокращающей точки (само)пересечения} и {\it
дефективного} класса Нильсена; (iii) {\it геометрическая точка
самопересечения}, {\it индекс} ее класса Нильсена (число Нильсена для
геометрических точек самопересечения оказывается более тонким
инвариантом, чем число Нильсена для точек само\-пе\-ре\-се\-чения, см.\ замечание~\ref {rem:mod2} и абзацы, расположенные до и после следствия~\ref {cor:geomselfint}) и естественно возникающее понятие {\it геометрически специальной} точки само\-пе\-ре\-се\-чения; (iv) {\it полуиндекс} класса Нильсена
геометрической точки самопересечения и необходимые для него понятия
{\it само\-со\-кра\-ща\-ющей} и {\it геомет\-ри\-чески
самосокращающей} точек само\-пе\-ре\-се\-чения, {\it дефективного
класса} геометрических точек само\-пе\-ре\-се\-чения. Мы показываем,
что соответствующие ``более тонкие'' числа Нильсена $NI^*[\gamma_1,\gamma_2]$ и
$NI^*[\gamma]$, к сожалению, совпадают с обычными (рассмотренными в \S\ref {sec:1}), поэтому даже для этих более тонких инвариантов аналог свойства Векена невыполнен.

Замкнутая кривая на многообразии называется {\it сохраняющей
ориентацию}, если локальная ориентация многообразия сохраняется при
непрерывном перенесении ориентации вдоль этой кривой. В противном
случае замкнутая кривая называется {\it меняющей ориен\-та\-цию}.

При общей деформации пары замкнутых кривых
$\gamma_1,\gamma_2\colon\,S^1\to S$ на поверхности $S$ точки
пересечения сохраняются, лишь слегка деформируясь, за исключением
некоторых значений параметра гомо\-то\-пии, при которых происходит
возникновение/уничтожение пары точек (см.\ рис.~1). При общей
деформации замкнутой кривой $\gamma\colon\,S^1\to S$ возможно также
возникновение/уничтожение одной точки самопересечения (см.\ рис.~2).

Окружность $S^1$ далее рассматривается как абелева группа $\RR/\ZZ$,
образ числа $t\in\RR$ при проекции $\RR\to\RR/\ZZ=S^1$ обозначается
через $\overline{t}$.

\begin{unDef}
\label{def:2.1}
 (a) Назовем {\it точкой пересечения} кри\-вых
$\gamma_1,\gamma_2\colon\,S^1\to S$ любую пару $(v_1,v_2)\in
S^1\times S^1$ со свойством $\gamma_1(v_1)=\gamma_2(v_2)$. Аналогично
{\it точкой самопересечения} кривой $\gamma\colon\,S^1\to S$
называется любая пара $(v_1,v_2)\in S^1\times S^1$, $v_1\ne v_2$, со
свойством $\gamma(v_1)=\gamma(v_2)$. Через $MI[\gamma_1,\gamma_2]$
обозначим минимальное число точек пересечения кривых $\gamma'_1$ и
$\gamma'_2$ среди всех пар кривых $(\gamma'_1,\gamma'_2)$, $\gamma'_i
\simeq \gamma_i$. Минимальное число $MI[\gamma]$ точек
самопересечения кривой $\gamma$ определяется как минимальное число
точек самопересечения среди всех кривых $\gamma'\simeq
\gamma$.\footnote{ В этом определении имеется своя тонкость, а
именно, для точки самопересечения $(v_1,v_2) \in S^1 \times S^1$
кривой $\gamma$ точка $(v_2,v_1)\in S^1 \times S^1$ тоже является
точкой самопересечения этой же кривой, и обе точки учитываются в
числе $MI[\gamma]$, хотя их образы в $S$ совпадают.}

\begin{figure}
\unitlength = 3.8mm 
\begin{center}
\begin{picture}(28,7.5)(-14,0)
 \thicklines
\qbezier[200](-15,0)(-12,3.5)(-15,7)
\qbezier[200](-11,0)(-14,3.5)(-11,7)

\put(0,.7) 
 {
\path(-10.5,3.5)(-8,3.5) \path(-10.5,3.5)(-10,3.8)
\path(-10.5,3.5)(-10,3.2) \path(-8,3.5)(-8.5,3.8)
\path(-8,3.5)(-8.5,3.2)
 }

\qbezier[200](-7,0)(-1.5,3.5)(-7,7)
  \put(-3.9,3.4){$\gamma\circ\alpha_2$}
\qbezier[200](-3,0)(-8.5,3.5)(-3,7)
  \put(-8.4,2.5){$\gamma\circ\alpha_1$}
  \put(-5.14,1.5){\put(-.14,-.07){$\bullet$}}
  \put(-4.5,1.6){\tiny$\gamma(v_1)=\gamma(v_2)$}
\put(-5.15,5.15){\put(-.14,-.07){$\bullet$}}
  \put(-4.5,5.2){\tiny$\gamma(v_1')=\gamma(v_2')$}
\path(-4.55,4.72)(-4.56,4.4) \path(-4.55,4.72)(-4.29,4.53)
\path(-5.38,4.72)(-5.37,4.4) \path(-5.38,4.72)(-5.64,4.53)
 \thinlines
\path(-5.5,3.05)(-4.45,3.95) \path(-5.2,2.4)(-4.5,3.05)
\path(-5.5,3.9)(-4.75,4.6)


\thicklines \qbezier[200](2,0)(2.5,0.5)(3,1.3)
\qbezier[200](6,0)(5.5,0.5)(5,1.3) \qbezier[200](3,1.3)(4,3)(5,1.3)

 \path(8,2)(6,2)
\path(8,2)(7.5,2.3) \path(8,2)(7.5,1.7) \path(6,2)(6.5,2.3)
\path(6,2)(6.5,1.7)

\qbezier[200](8,0)(10.1,1.5)(10.7,3)
\qbezier[200](10.7,3)(11.5,5.2)(10,5.5)
\put(11.3,3.7){$\gamma\circ\alpha$} \path(10,5.5)(10.42,5.58)
\path(10,5.5)(10.3,5.2) \qbezier[200](12,0)(9.9,1.5)(9.3,3)
\qbezier[200](9.3,3)(8.5,5.2)(10,5.5)
\put(9.85,1.6){\put(-.14,-.07){$\bullet$}}
\put(10.5,1.6){\tiny$\gamma(v)=\gamma(v')$}
 \thinlines
\path(9.3,4.1)(10.45,5.1) \path(9.5,3.35)(10.75,4.45)
\path(9.8,2.7)(10.6,3.45)
\end{picture}
\end{center}

\vskip 0.5truecm \noindent \!\!\!\!\!
\begin{tabular}{c} \small
Рис.~1. Возникновение эквивалентных\\ \small точек самопересечения
\end{tabular}
\begin{tabular}{c} \small
Рис.~2. Возникновение тривиальной\\ \small точки самопересечения
\end{tabular}
\end{figure}

\smallskip
(b)\quad Точки пересечения $(v_1,v_2)$ и $(v'_1,v'_2)$ кривых
$\gamma_1,\gamma_2$ называют {\it эквивалентными по Нильсену}, если
существуют пути $\alpha_i\colon\, [0,1] \to S^1$, $i = 1,2$, со
свойствами
 \begin{equation} \label {eq:equiv}
 \alpha_i(0) = v_i,\qquad \alpha_i(1) = v'_i, \qquad
 \gamma_1 \circ \alpha_1 \simeq_\partial  \gamma_2 \circ \alpha_2 .
 \end{equation}
Аналогично точки самопересечения $(v_1,v_2)$ и $(v'_1,v'_2)$ кривой
$\gamma$ назовем {\it эквивалентными по Нильсену}, если они являются
экви\-ва\-лент\-ными по Нильсену точками пересечения пары кривых
$\gamma,\gamma$. Если при этом $\alpha_1(t)\ne\alpha_2(t)$ для любого
$t\in[0,1]$, то точки самопересечения $(v_1,v_2)$ и $(v'_1,v'_2)$
назовем {\it строго экви\-ва\-лент\-ными по Нильсену}. Строго
эквивалентные точки самопересечения показаны на рис.~1 справа.

Точка самопересечения $(v_1,v_2)$ называется {\it тривиальной} или
{\it эквивалентной нулю}, если существует путь $\alpha\colon\, [0,1]
\to S^1$ со свойствами
 $$
\alpha(0) = v_1,\qquad \alpha(1) = v_2, \qquad \gamma \circ \alpha
\simeq_\partial 0 .
 $$
Тривиальная точка самопересечения показана на рис.~2 справа.
Тривиальные точки самопересечения образуют класс Нильсена, который
будем называть {\it тривиальным классом Нильсена}; он может быть
пустым.

Точку самопересечения $(v_1,v_2)$ кривой $\gamma$ назовем {\it
геометрически специальной}, если она эквивалентна
точке $(v_2,v_1)$. Если при этом соответствующая петля
$\gamma\circ\alpha_1$ сохраняет ориентацию для некоторой пары путей
$\alpha_1,\alpha_2$ в~(\ref {eq:equiv}), точку самопересечения
назовем {\it геометрически самосокращающей}.

\smallskip
(c)\quad Следуя~\cite[определение~2.2~(c,~d)]{BKZ}, мы назовем точку
пересечения $(v_1,v_2)$ кривых $\gamma_1,\gamma_2$ {\it
специальной}\footnote{ Определение специальной точки пересечения
допускает следующую переформулировку: точка пересечения является {\it
неспециальной} в том и только том случае, когда для любой
эквивалентной ей точки гомотопический класс (относительно концов)
каждого пути $\alpha_i$ на $S^1$ из определения~\ref {def:2.1}~(b)
определен однозначно, $i=1,2$, т.е.\ не зависит от выбора пары путей
$\alpha_1,\alpha_2$.}, если существует пара чисел $p,q\in\ZZ$,
$(p,q)\ne(0,0)$, такая, что замкнутые пути
$$
\delta_1(t):=\gamma_1(v_1+\overline{pt}) \qquad \mbox{и} \qquad
\delta_2(t):=\gamma_2(v_2+\overline{qt}), \qquad 0\le t\le 1,
$$
гомотопны относительно концов на поверхности $S$. Если при этом петля
$\delta_1$ меняет ориентацию для некоторой пары $p,q$, точка
пересечения называется {\it самосокращающей} (этот термин введен
в~\cite{Jez} для задачи о совпа\-де\-нии). Можно показать, что точка,
эквивалентная специальной или самосокращающей, тоже является
специальной или само\-со\-кра\-ща\-ющей соответственно;
поэтому кор\-ректно определены {\it специальные классы Нильсена} и
{\it дефективные классы Нильсена} точек пере\-се\-че\-ния (класс
Нильсена называется {\it дефективным}, если он состоит из
самосокращающих точек~\cite{Jez}).

Точку самопересечения кривой $\gamma$, а также ее класс Нильсена,
назовем {\it специальной},
если она является специальной точкой пере\-се\-чения пары кривых
$\gamma,\gamma$.
Ана\-ло\-гично точку самопересечения кривой $\gamma$ назовем {\it
самосокраща\-ю\-щей}, а ее класс Нильсена {\it дефективным}, если она
является самосокращающей точкой пересечения пары кривых
$\gamma,\gamma$.

Если в определении специальной (соответственно геометрически специальной)
точки самопересечения наложить условие $p=q\ne0$ (соответственно
условие строгой эквивалентности), получим понятия {\it строго
специ\-альной}, {\it строго самосокращающей} (соответственно {\it
строго гео\-мет\-ри\-чески специальной}, {\it строго геометрически
са\-мо\-со\-кра\-ща\-ю\-щей}) точек самопересечения. Можно показать,
что если точка само\-пе\-ре\-се\-че\-ния является строго специальной
(соответственно стро\-го самосокращающей, строго геометрически
специальной, строго геометрически самосокращающей), то все точки ее
класса Нильсена тоже являются таковыми; поэтому корректно определены
{\it строго специаль\-ные} классы Нильсена, 
{\it строго дефективные} классы Нильсена (т.е.\ состоящие из стро\-го самосокращающих точек) и {\it строго гео\-мет\-ри\-чески специальные} классы Нильсена точек самопересече\-ния.

\smallskip
(d)\quad {\it Индекс} изолированной точки и класса точек пересечения
кривых $\gamma_1,\gamma_2\colon\,S^1\to S$ на поверхности $S$
определим как в~(\ref {Gamma}) для $X_1=X_2=S^1$,
$Y=S$.\footnote{Можно показать, что индекс изолированной точки
пересечения кривых равен $0$, если в этой точке пересечение
несобственное, т.е.\ устранимо малым шевелением кривых. В случае
собственного пересечения индекс точки $(t_1,t_2)$ равен $1$, если
положительная полуветвь $\gamma_2|_{(t_2,t_2+\varepsilon)}$ второй
кривой лежит ``слева'' от ветви
$\gamma_1|_{(t_1-\varepsilon,t_1+\varepsilon)}$ первой кривой, и
равен $-1$ в противном случае. Здесь $\varepsilon$ -- малое
положительное число.} Индекс класса Нильсена точек пересечения
принимает значения или в группе $\ZZ$ в случае ориентированной
поверхности $S$, или в группе $\ZZ_2=\ZZ/2\ZZ$ в случае
неориентируемой поверхности $S$ (а также в случае {\it индекса по
модулю 2}), или в $\ZZ_+=\{k\in\ZZ\mid k\ge0\}$ в случае так
называемого {\it полуиндекса}~\cite{DobJez,Jez}, определенного чуть
ниже. Класс Нильсена точек пересечения называется {\it существенным},
если его индекс не нуль; {\it число Нильсена} $NI[\gamma_1,\gamma_2]$
точек пересечения является числом существенных классов. Число классов
Нильсена, имеющих ненулевой полуиндекс, обозначим через
$NI^*[\gamma_1,\gamma_2]$. Ясно, что
 $$
 MI[\gamma_1,\gamma_2]\ge NI^*[\gamma_1,\gamma_2]\ge NI[\gamma_1,\gamma_2].
 $$
Если числа $MI[\gamma_1,\gamma_2]$ и $NI[\gamma_1,\gamma_2]$
совпадают, мы скажем, что выполняется {\it свойство Векена для задачи
о пересечении} (см.\ замечание~\ref {rem:mod2} ниже).

Если класс Нильсена дефективен (см.~(c)), его {\it полуиндекс}
определяется как образ $\mod2$-индекса при вложении
$\ZZ_2=\ZZ/2\ZZ\hookrightarrow\ZZ_+$, переводящем $2\ZZ\mapsto0$,
$1+2\ZZ\mapsto1$. Для недефективного класса Нильсена точек
пересечения кривых $\gamma_1,\gamma_2$, состоящего из изолированных
точек пересечения, {\it полуиндекс} определяется как абсолютное
значение суммы индексов точек пересечения из этого класса, где
индексы точек пересечения определяются при помощи ``согласованных''
локальных ориентаций поверхности. А именно, локальная ориентация в
точке $\gamma_1(v_1')=\gamma_2(v_2')$, с помощью которой определяется
индекс точки пересечения $(v_1',v_2')$, получается из локальной
ориентации в точке $\gamma_1(v_1)=\gamma_2(v_2)$, с помощью которой
определяется индекс точки пересечения $(v_1,v_2)$, перенесением вдоль
соответствующего пути $\gamma_1\circ\alpha_1$, см.~(\ref {eq:equiv}).
Полуиндекс класса Нильсена не зависит от выбора пары путей $\alpha_i$
из~(b) в силу недефективности этого класса.

Для замкнутой кривой $\gamma$ на $S$ таким же путем определяются
понятия {\it индекса} нетривиального (см.~(b)) класса Нильсена точек
само\-пе\-ре\-се\-чения (индекс тривиального класса полагается
равным нулю), {\it существенного класса} точек самопересечения, {\it
числа Ниль\-сена} $NI[\gamma]$ точек самопересечения, {\it свойства
Векена для задачи о самопересечении}. {\it Полуиндекс} нетривиального
класса Нильсена точек само\-пе\-ре\-се\-чения определяется с
использова\-нием понятия дефективного класса точек
само\-пе\-ре\-се\-чения. При рассмотрении классов строгой
эквивалентности по Нильсену точек самопересечения (см.\ (b))
получаются аналогичные понятия, при этом для определения {\it
полу\-ин\-декса} класса строгой эквивалент\-ности используется
понятие строго дефективного класса Нильсена точек
само\-пе\-ре\-се\-чения, см.~(c). Число классов строгой
эквивалентности по Нильсену, имеющих ненулевой полуиндекс, обозначим
через $NI^*[\gamma]$. Конечно,
 $$
 MI[\gamma]\ge NI^*[\gamma]\ge NI[\gamma]
 $$
и числа Нильсена $NI[\gamma]$ и $NI^*[\gamma]$ являются инвариантами
свобод\-ного гомотопического класса замкнутой кривой.

\smallskip
(e)\quad Пусть $(v_1,v_2)$ --- точка самопересечения кривой $\gamma$.
Неупо\-ря\-до\-ченную пару $\{v_1,v_2\}$ назовем {\it геометрической
точкой само\-пе\-ре\-се\-чения} $\gamma$. Геометрические точки
самопересечения $\{v_1,v_2\}$ и $\{v_1',v_2'\}$ назовем {\it
эквивалентными по Нильсену}, если точка $(v_1,v_2)$ эквивалентна
$(v_1',v_2')$ или $(v_2',v_1')$. {\it Индекс} геометрической точки
самопересечения и ее класса Нильсена определим по модулю 2 (т.е.\
индекс геометрической точки равен $1$, если самопересечение
собственное, и индекс равен $0$, если самопересечение несобственное).
Обозначим через $MI_\geom[\gamma]$ минимальное число геометрических
точек самопересечения кривой $\gamma'$ среди всех кривых
$\gamma'\simeq \gamma$ (очевидно,
$MI_\geom[\gamma]=\frac12MI[\gamma]$); через $NI_\geom[\gamma]$ --
число сущест\-венных классов Нильсена геометрических точек
самопересечения $\gamma$.

{\it Полуиндекс} класса Нильсена геометрических точек
само\-пе\-ре\-се\-чения определим аналогично (d). А именно, класс
Нильсена геометрической точки самопересечения $\{v_1,v_2\}$ назовем
{\it дефек\-тив\-ным}, если точка самопересечения $(v_1,v_2)$
является само\-со\-кра\-ща\-ющей или геометрически самосокращающей
(см.~(b,~c)). Для дефективных классов определим полуиндекс при помощи
$\mod2$-индекса, а для недефективных классов определим полуиндекс
так. Полуиндекс класса Нильсена геометрической точки самопересечения
$\{v_1,v_2\}$ определяется либо как полуиндекс класса Нильсена точки
самопересечения $(v_1,v_2)$ (см.\ (d)), если точка $(v_1,v_2)$ не является
гео\-мет\-ри\-чески специальной (см.\ (b)), либо как половина этого
полуиндекса, если точка $(v_1,v_2)$ является геометрически
специальной.

Аналогично определяются классы {\it строго эквивалентных по Нильсену}
геометрических точек самопересечения; {\it индекс} такого класса
определяется по модулю 2, а {\it полуиндекс} определяется так.
Класс Нильсена геометрической точки самопе\-ре\-сечения $\{v_1,v_2\}$ назовем {\it строго
дефективным}, если точка самопересечения $(v_1,v_2)$ является строго
само\-со\-кра\-ща\-ющей или строго геометрически самосокращающей,
см.~(c). 
Если класс строгой эквивалентности геометрической точки самопересечения $\{v_1,v_2\}$ строго дефективен, его {\it полуиндекс} определяется при помощи $\mod2$-индекса. 
Если этот класс 
не является строго дефективным, его {\it полуиндекс} определяется либо как полуиндекс класса строгой эквивалентности точки самопересечения $(v_1,v_2)$ (см.\ (d)), если точка $(v_1,v_2)$ не является строго гео\-мет\-ри\-чески специальной (см.\ (c)), либо как половина этого
полуиндекса, если точка $(v_1,v_2)$ является строго геометрически
специальной. Число классов строгой эквивалентности геометрических точек,
имеющих ненулевой полуиндекс, обозначим через $NI^*_\geom[\gamma]$. Конечно,
 $$
 MI_\geom[\gamma]\ge NI_\geom^*[\gamma]\ge NI_\geom[\gamma]
 $$
и геометрические числа Нильсена $NI_\geom[\gamma]$ и
$NI_\geom^*[\gamma]$ явля\-ют\-ся инвариантами свободного
гомотопического класса замкнутой кривой.
\end{unDef}

Определения~\ref {def:2.1}~(b,~c) {\it эквивалентных} точек
(само)пересечения, {\it специальной} и {\it самосокращающей} точек
(само)пересечения, а также {\it тривиальной}, {\it геометрически
специальной} и {\it геометри\-чески самосокращающей} точек
самопересечения и их ``стро\-гих'' аналогов допускают следу\-ю\-щую
пере\-фор\-му\-лировку.

\begin{unLem}\label{lem:topol}
Пусть $\gamma,\gamma_1,\gamma_2\colon\,[0,1]\to S$ --- замкнутые
кривые.

{\rm(a)} Каждой точке пересечения $(v_1,v_2)$ кривых
$\gamma_1$ и $\gamma_2$, $v_1,v_2\in[0,1)$, сопо\-ста\-вим путь
$\gamma_{v_1,v_2}=\gamma_1|_{[0,v_1]}\gamma_2|_{[0,v_2]}^{-1}$,
ведущий из точки $\gamma_1(0)$ в точку $\gamma_2(0)$ на $S$. Тогда
верны следующие критерии: {\rm
$$
 (v_1,v_2),(v_1',v_2')\mbox{ эквивалентны} \iff \exists p,q\in\ZZ :
 \gamma_{v_1,v_2} \simeq_\partial \gamma_1^p\gamma_{v_1',v_2'} \gamma_2^q ;
$$
$$
 (v_1,v_2)\mbox{ специальна} \iff \exists (p,q)\in\ZZ^2\setminus\{(0,0)\} :
 \gamma_{v_1,v_2} \simeq_\partial \gamma_1^p\gamma_{v_1,v_2} \gamma_2^q;
$$
$$
 (v_1,v_2)\mbox{ самосокр.} \iff (v_1,v_2)\mbox{ специальна и $\gamma_1^p$ меняет ориентацию}.
$$
}%

{\rm(b)} Каждой точке само\-пе\-ре\-се\-че\-ния $(v_1,v_2)$ кривой
$\gamma$, $v_1,v_2\in[0,1)$, сопо\-ста\-вим замкнутый путь
$\gamma_{v_1,v_2}:=\gamma|_{[0,v_1]}\gamma|_{[0,v_2]}^{-1}$ на $S$.
Тогда $\gamma_{v_1,v_2}\simeq_\partial\gamma_{v_2,v_1}^{-1}$ и верны
следующие критерии: {\rm
$$
 (v_1,v_2),(v_1',v_2')\mbox{ эквивалентны} \iff \exists p,q\in\ZZ :
 \gamma_{v_1,v_2} \simeq_\partial \gamma^p\gamma_{v_1',v_2'} \gamma^q ;
$$
$$
 (v_1,v_2)\mbox{ тривиальна} \iff
 \gamma_{v_1,v_2} \simeq_\partial \gamma^p \mbox{ для некоторого } p\in\ZZ;
$$
$$
 (v_1,v_2)\mbox{ специальна} \iff \exists (p,q)\in\ZZ^2\setminus\{(0,0)\} :
 \gamma_{v_1,v_2} \simeq_\partial \gamma^p\gamma_{v_1,v_2} \gamma^q;
$$
$$
 (v_1,v_2)\mbox{ геом.\ специальна} \iff \exists p,q\in\ZZ :
 \gamma_{v_1,v_2} \simeq_\partial \gamma^p\gamma_{v_1,v_2}^{-1} \gamma^q ;
$$
$$
 (v_1,v_2)\mbox{ самосокр.} \iff (v_1,v_2)\mbox{ специальна и $\gamma^q$ меняет ориентацию};
$$
$$
 (v_1,v_2)\mbox{ геом.\ самосокр.} \iff
 \left\{\begin{array}{l}(v_1,v_2)\mbox{ геом.\ специальна и}\\ \gamma^{-q}\gamma_{v_1,v_2} \mbox{ сохраняет ориентацию;}
 \end{array}
 \right.
$$
$$
 (v_1,v_2),(v_1',v_2')\mbox{ эквивалентны$^*$} \iff \exists q\in\ZZ :
 \gamma_{v_1,v_2} \simeq_\partial \gamma^{-\hat q}\gamma_{v_1',v_2'} \gamma^q ;
$$
$$
 (v_1,v_2)\mbox{ специальна$^*$} \iff \exists q\in\ZZ\setminus\{0\} :
 \gamma_{v_1,v_2} \simeq_\partial \gamma^{-q}\gamma_{v_1,v_2} \gamma^q;
$$
$$
 (v_1,v_2)\mbox{ геом.\ специальна$^*$} \iff \exists q\in\ZZ :
 \gamma_{v_1,v_2} \simeq_\partial \gamma^{\eta-q}\gamma_{v_1,v_2}^{-1} \gamma^q ;
$$
$$
 (v_1,v_2)\mbox{ самосокр.$^*$} \iff (v_1,v_2)\mbox{ специальна$^*$ и $\gamma^q$ меняет ориент.};
$$
$$
 (v_1,v_2)\mbox{ геом.\ самосокр.$^*$} \iff
 \left\{\begin{array}{l}(v_1,v_2)\mbox{ геом.\ специальна$^*$ и}\\ \gamma^{-q}\gamma_{v_1,v_2} \mbox{ сохраняет ориентацию,}
 \end{array}
 \right.
$$
}%
где последние пять критериев относятся к строго экви\-ва\-лент\-ным
по Нильсену точкам самопересечения и соответст\-ву\-ющим понятиям:
строго специальной, строго геометрически специальной и т.п.\ точкам
самопересечения, где звездочка $*$ означает ``строго''. Здесь
$\eta:=\sgn(v_2-v_1)$, $\eta_i:=\sgn(v_i'-v_i+\varepsilon)$, $i=1,2$, $0<\varepsilon\ll1$,
$\hat q:=q+(\eta_1-\eta_2)/2-1$ в случае рас\-по\-ло\-же\-ния точек на отрезке $[0,1]$
в циклическом порядке $v_1,v_1',v_2',v_2$ (возможно $v_1=v_1'$),
$\hat q:=q+(\eta_1-\eta_2)/2+1$ при расположении точек в обратном циклическом порядке (возможно $v_2=v_2'$), $\hat q:=q+(\eta_1-\eta_2)/2$ в остальных случаях.
\qed
\end{unLem}

В настоящей работе нам понядобятся лишь критерии эквивалентности точек (само)пересечения (см.\ первую формулу леммы~\ref {lem:topol}~(a,~b)), доказательство которых не представляет труда. Переформулировки остальных понятий приведены в лемме~\ref {lem:topol} для полноты изложения.

\begin{nnExa}
Из леммы~\ref {lem:topol} легко видно, что тривиальная точка
самопересечения кривой $\gamma$ является специальной, строго
спе\-ци\-аль\-ной, геометрически специ\-альной и геометрически
самосокра\-ща\-ющей.
Тривиальная точка является строго геометрически специаль\-ной
(соответственно строго геометрически самосокращающей,
само\-сокращающей, строго самосокращающей) тогда и только тогда,
когда кривая $\gamma$ стягиваема (соответственно стягиваема, меняет
ориентацию, меняет ориентацию).
\end{nnExa}

\begin{unRem}\label{rem:invol}
{\sl Инволюция $(v_1,v_2)\mapsto(v_2,v_1)$ на множестве 
$\int(\gamma)$ точек самопересечения} пере\-во\-дит классы Нильсена в
классы Нильсена, поэтому она {\sl инду\-ци\-рует инволюцию на
множестве классов Нильсена} точек самопересечения. Эта инволюция
переводит сущест\-венные классы Нильсена в сущест\-венные,
специ\-аль\-ные в специаль\-ные, дефективные в дефективные, а
тривиальный класс оставляет на месте. Заметим, что любой класс,
оставляемый инволюцией на месте, является $\mod2$-несуществен\-ным
(так как точки в таком классе возникают парами:
$(v_1,v_2)\sim(v_2,v_1)$). Поэтому геометрически специальные (т.е.\
инвариантные относительно инво\-лю\-ции
$\int(\gamma)\to\int(\gamma)$, $(v_1,v_2)\mapsto(v_2,v_1)$) классы
Нильсена точек самопересечения не являются $\mod2$-существен\-ными.
Отметим, что тривиальный класс по нашему определению несуществен.
\end{unRem}

\begin{unRem} \label {rem:mod2}
Если отличные от нуля полуиндексы классов Нильсена равны $1$ (что, в
действительности, имеет место по теоремам \ref {pro:2.4} и \ref
{pro:2.2}), то понятие существенного класса Нильсена (а потому и
число Нильсена) совпадает для рассмотренных в определении~\ref
{def:2.1}~(d) трех типов индекса. Для чисел Нильсена, определенных с
помощью индексов в группе $\ZZ_2$, легко проверяется цепочка
неравенств
\begin{equation} \label{geom:Wecken}
MI[\gamma] = 2 MI_\geom[\gamma] \ge 2 NI_\geom[\gamma] \ge NI[\gamma]
\end{equation}
(последнее неравенство следует из того, что $\mod2$-существенные
классы точек самопересечения не инвариантны относительно
инво\-лю\-ции $(v_1,v_2)\mapsto(v_2,v_1)$, см.\ замечание~\ref
{rem:invol}). Значит, для любой замкнутой кривой $\gamma$,
обладающей свойством Векена $MI[\gamma]=NI[\gamma]$ для точек
самопересечения, верно равенство $NI[\gamma]=2NI_\geom[\gamma]$ и
выполнено {\it свойство Векена для геометрических точек
самопересечения}: $MI_\geom[\gamma]=NI_\geom[\gamma]$.
\end{unRem}

\begin{unDef}\label{defrem:Reidem}
(a) {\it Классами Райдемайстера} замкнутой кривой $\gamma$ назовем
двусторонние смежные классы $Z_\gamma gZ_\gamma$,
$g\in\pi_1(S,\gamma(0))$, группы $\pi_1(S,\gamma(0))$ по циклической
подгруппе $Z_\gamma$, порож\-ден\-ной гомотопическим классом кривой
$\gamma$. Число $RI[\gamma]$ классов Райдемайстера назовем {\it
числом Райдемайстера для самопересе\-чения} замкнутой кривой
$\gamma$. Подмножества вида $(Z_\gamma gZ_\gamma)\cup(Z_\gamma
g^{-1}Z_\gamma)$ группы $\pi_1(S,\gamma(0))$ назовем {\it
геометрическими классами Райде\-май\-стера для самопересечения}
замкнутой кривой $\gamma$. Число геомет\-ри\-ческих классов
Райдемайстера обозначим через $RI_\geom[\gamma]$.

(b) Для задачи о пересечении можно ввести похожие понятия: {\it
Классы Райдемайстера пары замкнутых кривых} $\gamma_1,\gamma_2$ с
$\gamma_1(0)=\gamma_2(0)$ являются  двусторонними смежными классами
$Z_{\gamma_1} gZ_{\gamma_2}$ по циклическим подгруппам
$Z_{\gamma_1},Z_{\gamma_2}$. Число $RI[\gamma_1,\gamma_2]$ классов
Райде\-май\-стера назовем {\it числом Райдемайстера для пересечения}
замкнутых кривых $\gamma_1,\gamma_2$.
 \end{unDef}

Из первой формулы леммы~\ref {lem:topol}~(a,~b) следует, что для
любой пары кривых $\gamma_1',\gamma_2'$, такой, что
$\gamma_1'\simeq\gamma_1$, $\gamma_2'\simeq\gamma_2$ и
$\gamma_1'(0)=\gamma_2'(0)$, имеется взаимно-однозначное соответствие
между классами Нильсена точек пересечения $\gamma_1',\gamma_2'$ и
некоторыми классами Райдемайстера пары кривых $\gamma_1,\gamma_2$.
Аналогично, для любой кривой $\gamma'\simeq\gamma$ имеется
взаимно-однозначное соответствие между классами Нильсена точек
самопересечения $\gamma'$ и некоторыми классами Райдемайстера кривой
$\gamma$, и аналогично для геометрических классов. Следовательно,
   $$
NI[\gamma]\le RI[\gamma], \quad NI[\gamma_1,\gamma_2]\le
RI[\gamma_1,\gamma_2], \quad NI_\geom[\gamma]\le RI_\geom[\gamma].
   $$

\section{ЗАМКНУТЫЕ КРИВЫЕ НА ЦИЛИНДРЕ И ЛИСТЕ МЕБИУСА}\label{sec:3}

Для $k\in\ZZ$ положим
 \begin{equation} \label{eq:delta}
\delta(k) := \frac{(-1)^k+1}2=\left\{ \begin{array}{ll}
 0, & \mbox{если $k$ нечетно}, \\
 1, & \mbox{если $k$ четно}.
\end{array} \right.
 \end{equation}

Следующие леммы показывают, что свойство Векена может нарушаться в
задачах о точках пересечения и самопересечения кривых на цилиндре и
листе Мебиуса. Эти леммы являются ключевыми в доказательстве
теорем~\ref {pro:2.4}~(b) и~\ref {pro:2.2}~(b).

\begin{unLem}\label{lem:Cylinder}
Пусть $\gamma,\gamma_1,\gamma_2\colon\, S^1 \to C$ -- замкнутые
кривые на цилиндре
 $C=S^1\times[-1,1]$, и пусть $k=\deg(\rho\circ\gamma)\ne0$,
 $k_i=\deg(\rho\circ\gamma_i)$, где $\rho\colon\,C\to S^1$ -- проекция. Тогда
любой класс Нильсена точек самопересечения или пересечения
несуществен (и даже недефективен и имеет полуиндекс $0$), и
   $$
 MI[\gamma] = 2|k|-2, \qquad NI^*[\gamma] = NI[\gamma] = 0, \qquad RI[\gamma] = |k|;
   $$
   $$
 MI[\gamma_1,\gamma_2] = NI^*[\gamma_1,\gamma_2] = NI[\gamma_1,\gamma_2] = 0,
   $$
   $$
 RI[\gamma_1,\gamma_2] = \left\{ \begin{array}{ll}
                               \!\!\!\infty,        & \mbox{если } k_1=k_2=0, \\
                               \!\!\!\gcd(k_1,k_2), & \mbox{если }|k_1|+|k_2|>0.
                                \end{array} \right.
   $$
Более того, кривая $\gamma$ имеет не менее $|k|-1$ классов Нильсена
точек самопересечения, каждый из которых не совпадает с тривиальным и
содержит не менее двух (строго эквива\-лент\-ных, {\rm см.\
определение~\ref {def:2.1}~(b)}) точек самопересечения, хотя
несуществен (и даже недефективен и имеет полуин\-декс $0$); в
точности $\delta(k)$ из этих классов является геометри\-чески
специальным (и даже состоящим из строго геометри\-чески
само\-со\-кра\-ща\-ю\-щих точек), {\rm см.\ определение~\ref
{def:2.1}~(b,~c) и (\ref {eq:delta})}.
 \end{unLem}

\begin{unCor}\label{cor:geomCylinder}
В условиях леммы~$\ref {lem:Cylinder}$ для замкнутой кривой
$\gamma\colon\, S^1 \to C$ на цилиндре отличные от нуля полуиндексы
классов геометрических точек самопересечения равны $1$ и
   $$
MI_\geom[\gamma] = |k|-1, \qquad NI_\geom^*[\gamma] =
NI_\geom[\gamma] = \delta(k),
 $$
 $$
RI_\geom[\gamma] = [\frac{|k|}{2}]+1.
   $$
Более того, кривая $\gamma$ имеет не менее $[\frac{|k|}{2}]$ классов
Нильсена гео\-мет\-ри\-чес\-ких точек самопересечения, отличных от
триви\-аль\-ного; в точности $\delta(k)$ из этих классов является
(строго) гео\-мет\-ри\-чески специальным (он существен
и строго дефе\-к\-ти\-вен), а каждый из остальных $[\frac{|k|-1}2]$
классов содержит не ме\-нее двух (строго эквивалентных)
геометрических точек, хотя не\-су\-щест\-вен (и даже недефективен
и имеет полуиндекс $0$). \qed
\end{unCor}

\PROOF леммы \ref {lem:Cylinder}. Так как разные кривые можно
``развести'' на разные компоненты края цилиндра,
$MI[\gamma_1,\gamma_2]=0$. Следовательно, $NI[\gamma_1,\gamma_2]=0$.
Так как группа $\pi_1(C)=\ZZ$ абелева, число Райдемайстера
$RI[\gamma_1,\gamma_2]$ равно числу элементов факторгруппы группы
$\ZZ$ по подгруппе, порожденной числами $k_1,k_2$. Отсюда следует,
что $RI[\gamma_1,\gamma_2]=\gcd(k_1,k_2)$ при $(k_1,k_2)\ne(0,0)$,
$RI[\gamma_1,\gamma_2]=\infty$ для $k_1=k_2=0$. Решение задачи о
самопересечении разобьем на несколько шагов.

{\it Шаг 1.} Рассмотрим какую-нибудь кривую $\gamma'\simeq\gamma$,
имеющую $MI[\gamma]$ точек самопересечения, т.е.\ минимальное число
точек само\-пе\-ре\-се\-че\-ния. Каждую кривую $\gamma'$ на цилиндре,
имеющую конечное число точек самопересечения, можно продеформировать
--- с сохранением точек самопересечения и их разбиения на классы
Нильсена
--- так, чтобы точки $\gamma'(0)$ и $\gamma'(1/2)$ принадлежали
разным компонентам края цилиндра. Поэтому мы можем и будем считать,
что точки $\gamma'(0)$ и $\gamma'(1/2)$ принадлежат разным
компонентам края цилиндра. Кривая $\gamma'$ состоит из двух
незамкнутых дуг $\delta_1:=\gamma'|_{[0,1/2]}$ и
$\delta_2:=\gamma'|_{[1/2,1]}$, где начало одной дуги является концом
другой и наоборот, и эти точки принадлежат разным компонентам края
цилиндра.

Для вычисления минимального числа $MI[\gamma]$ точек
само\-пе\-ре\-се\-чения кривой $\gamma$ рассмотрим задачу о
минимальном числе $MI[\delta_1,\delta_2]$ точек пересечения дуг
$\delta_1',\delta_2'$ в гомотопических классах дуг
$\delta_1,\delta_2$ с фиксированными концами; в процессе гомотопии
концы дуг остаются неподвижными и, в частности, принадлежащими краю
цилиндра. Заметим, что $\delta_1\not\simeq_\partial\delta_2$, так как
$k\ne0$. Известно, что в этой задаче о пересечении справедливо
свойство Векена, а мини\-маль\-ное число точек пересечения в данных
классах гомотопии имеют кривые, являющиеся кратчайшими относительно
некоторой плоской римановой метрики на цилиндре, см.,
например,~\cite[теорема~II]{Turaev} или~\cite[theorem~2.11]{GKZ2},
где $MI[\delta_1,\delta_2]=k(\delta_1,\delta_2)$,
$NI[\delta_1,\delta_2]=
|\lambda(\delta_1,\delta_2)|=|\tilde\lambda(\delta_1,\delta_2)|$.
Отсюда легко находим
 $$
MI[\delta_1,\delta_2]=NI[\delta_1,\delta_2]=|k|+1.
 $$
Следовательно, дуги $\delta_1$ и $\delta_2$ имеют не менее $|k|+1$
точек пересечения, из которых две точки $(0,1)$ и $(1/2,1/2)$
отвечают концам дуг. Так как любой точке пересечения кроме указанных
двух точек отвечает пара точек самопересечения кривой $\gamma'$, то
$MI[\gamma]=2(MI[\delta_1,\delta_2]-2)$. Отсюда следует требуемое
равенство $MI[\gamma]=2(|k|-1)$.

{\it Шаг 2.} Для вычисления числа Нильсена $NI[\gamma]$ сопоставим
каждой точке самопересечения $(v_1,v_2)$ кривой $\gamma$ замкнутый
путь $\gamma_{v_1,v_2}=\gamma|_{[0,v_1]}\cdot\gamma|_{[0,v_2]}^{-1}$
и целое число $g_{v_1,v_2}:=\deg(\rho\circ\gamma_{v_1,v_2})$, где
$\rho\colon\,C\to S^1$ -- проекция. С учетом первой формулы
леммы~\ref{lem:topol}~(b) и изоморфизма
$\rho_\#\colon\,\pi_1(C)\to\pi_1(S^1)\cong\ZZ$, точки самопересечения
$(v_1,v_2)$ и $(v_1',v_2')$ эквивалентны в том и только том случае,
когда $k\mid(g_{v_1,v_2}-g_{v_1',v_2'})$. Следовательно, классы
Нильсена точек самопересечения $(v_1,v_2)$ кривой $\gamma$ находятся
во взаимно-однозначном соответствии с образами чисел $g_{v_1,v_2}$
при проекции $\ZZ\to\ZZ/|k|\ZZ$, $i\mapsto i\mod|k|$.

Рассмотрим кратчайшие дуги $\delta_1',\delta_2'$ в гомотопических
классах
   $$
\delta_1'\simeq_\partial\delta_1,\qquad
\delta_2'\simeq_\partial\delta_2
   $$
относительно какой-нибудь плоской римановой метрики на цилиндре.
Соответствующая замкнутая кривая $\gamma'=\delta_1'\cdot\delta_2'$
имеет $2(|k|-1)$ точек самопересечения; отвечающие им $2(|k|-1)$
чисел $g_{v_1,v_2}$ образуют множество
$\{\pm1,\pm2,\dots,\pm(|k|-1)\}$; при этом точка, отвечающая числу
$\pm i$, имеет индекс $\pm1$ при подходящем выборе ориентации
цилиндра. Так как указанное множество чисел распадается на пары чисел
вида $\{i,-(|k|-i)\}$, равных в группе $\ZZ/|k|\ZZ$, то отвечающие
каждой паре чисел точки самопересечения кривой $\gamma'$ образуют
класс Нильсена и имеют противоположные индексы. Поэтому все классы
Нильсена несущественны, т.е.\ $NI[\gamma]=0$. Так как цилиндр
ориентируем, все классы недефективны, а потому имеют полуиндекс $0$.

{\it Шаг 3.} Покажем, что для любого $i\in\{1,2,\dots,|k|-1\}$ кривая
$\gamma$ имеет не менее двух точек самопересечения в классе Нильсена,
отвечающем элементу $i\mod|k|\in\ZZ/|k|\ZZ$, см.\ шаг~2. Как на
шаге~1, каждую (но уже не обязательно имеющую $MI[\gamma]$ точек
самопересечения) кривую $\gamma''\simeq\gamma$, имеющую конечное
число точек самопересечения, можно продеформировать --- с сохранением
точек самопересечения и их разбиения на классы Нильсена --- так,
чтобы точки $\gamma''(0)$ и $\gamma''(1/2)$ принадлежали разным
компонентам $\partial C$, после чего можно разбить ее на две (не
обязательно простые) дуги: $\gamma''=\delta_1''\cdot\delta_2''$.
Каждая точка самопересечения кривой $\gamma''$ отвечает либо точке
самопересечения одной из дуг $\delta_1''$ или $\delta_2''$, либо
точке пересечения ровно одной из пар дуг $(\delta_1'',\delta_2'')$
или $(\delta_2'',\delta_1'')$; эквивалентным точкам пересечения любой
из пар дуг отвечают эквивалентные (и даже строго эквивалентные) точки
самопересечения кривой $\gamma''$.

Аналогично шагу 2 любой точке пересечения $(v_1,v_2)$ пары дуг
$\delta_1'',\delta_2''$ сопоставим замкнутый путь
$\delta_{v_1,v_2}=\delta_1''|_{[0,v_1]}\cdot\delta_2''|_{[1/2,v_2]}^{-1}$
и целое число $g_{v_1,v_2}:=\deg(\rho\circ\delta_{v_1,v_2})$.
Известно~\cite{GKZ2}, что точки пересечения $(v_1,v_2)$ и
$(v_1',v_2')$ дуг $\delta_1'',\delta_2''$ эквивалентны в том и только
том случае, когда $\delta_{v_1,v_2}\simeq_\partial\delta_{v_1',v_2'}$
(последнее равносильно равенству $g_{v_1,v_2}=g_{v_1',v_2'}$);
поэтому указанное сопоставление индуцирует инъективное сопоставление
целых чисел классам Нильсена точек пересечения пары дуг
$\delta_1'',\delta_2''$; класс эквивалентности точки пересечения
$(v_1,v_2)$ будем называть {\it $g_{v_1,v_2}$-м классом}.

Рассмотрим описанное выше сопоставление целых чисел классам Нильсена
точек пересечения пары дуг $\delta_1'',\delta_2''$ (соответственно
пары дуг $\delta_2'',\delta_1''$). Тогда существенным классам
сопоставляются числа $1,2,\dots,|k|-1$ (соответственно числа
$-1,-2,\dots,-(|k|-1)$). При этом точке пересечения $i$-го класса
пары $(\delta_1'',\delta_2'')$ и точке пересечения $-(|k|-i)$-го
класса пары $(\delta_2'',\delta_1'')$ отвечают две разные точки
самопересечения кривой $\gamma''$ одного класса Нильсена $i\mod|k|$
($1\le i\le|k|-1$). Непосредственная проверка показывает, что такие
две точки самопересечения кривой $\gamma''$ строго эквивалентны. Так
как $i$ -- любое целое число в интервале $1\le i\le|k|-1$, то кривая
$\gamma''$ имеет не менее $|k|-1$ классов Нильсена точек
самопересечения, причем в каждом классе не менее двух строго
эквивалентных точек.

{\it Шаг 4.} В силу второй формулы леммы \ref {lem:topol}~(b) точки
самопере\-се\-че\-ния, отвечающие рассматриваемым классам,
нетривиальны, так как $k\nmid g_{v_1,v_2}$. Точка самопересечения
$(v_1,v_2)$ геометрически специальна (а потому и геометрически
самосокращающая ввиду ориентируемости цилиндра), если
$g_{v_1,v_2}=|k|/2$. Последнее выполнено при четном $k$ для
$|k|/2$-го класса пары $(\delta_1'',\delta_2'')$ и для $-|k|/2$-го
класса пары $(\delta_2'',\delta_1'')$, а потому при четном $k$
существует геометрически специальный класс точек самопересечения
кривой $\gamma''$. Любая его точка $(v_1,v_2)$ строго геометрически
специальна (а потому и строго геометрически самосокращающая ввиду
ориентируемости цилиндра), так как она строго эквивалентна точке
$(v_2,v_1)$, см.\ выше.

Так как группа $\pi_1(C)=\ZZ$ абелева, число Райдемайстера
$RI[\gamma]$ равно числу элементов факторгруппы группы $\ZZ$ по
подгруппе, порожденной числом $k$. Отсюда следует, что
$RI[\gamma]=|k|$. \qed \ENDPROOF

\begin{unLem}\label{lem:Moebius}
Пусть $\gamma,\gamma_1,\gamma_2\colon\, S^1 \to M$ -- замкнутые
кривые на листе Мебиуса $M$, и пусть $k=\deg(\rho\circ\gamma)\ne0$,
$k_i=\deg(\rho\circ\gamma_i)$, где $\rho\colon\,M\to S^1$ --
стандартное расслоение со слоем отрезок. Тогда \ $RI[\gamma] = |k|$,
   $$
MI[\gamma] = \left\{
\begin{array}{ll} |k|-2, & 2\mid k, \\
                  |k|-1, & 2\nmid k;
\end{array} \right.
\qquad NI^*[\gamma] = NI[\gamma] = \left\{
\begin{array}{ll}     0, & 2\mid k, \\
                  |k|-1, & 2\nmid k;
\end{array} \right.
   $$
   $$
\begin{array}{rl}
MI[\gamma_1,\gamma_2] =& \left\{
\begin{array}{ll}     0, & \mbox{если } 2\mid k_1k_2, \\
    \min\{|k_1|,|k_2|\}, & \mbox{если }2\nmid k_1k_2;
\end{array} \right.
 \\
NI^*[\gamma_1,\gamma_2] = NI[\gamma_1,\gamma_2] =& \left\{
\begin{array}{ll}     0, & \mbox{если }2\mid k_1k_2, \\
          \gcd(k_1,k_2), & \mbox{если }2\nmid k_1k_2;
\end{array} \right.
 \\
RI[\gamma_1,\gamma_2] =& \left\{
\begin{array}{ll}
    \infty,        & \mbox{если }k_1=k_2=0, \\
    \gcd(k_1,k_2), & \mbox{если }|k_1|+|k_2|>0.
\end{array} \right.
\end{array}
   $$

Если $k$ нечетно, то выполнено свойство Векена, все классы Нильсена
точек самопересечения строго дефективны и не являются геометрически
специальными, а все классы кроме тривиального существенны {\rm (см.\
определение~\ref {def:2.1}(b,~c))}.

Если $k$ четно, то любой класс точек самопересечения несуществен (и
даже недефективен и имеет полуиндекс $0$); су\-щест\-вует не менее
$|\frac{k}{2}|-1$ классов Нильсена, каждый из кото\-рых не совпадает
с тривиальным и содержит не менее двух (строго экви\-ва\-лент\-ных)
точек самопересечения, хотя несуществен (и даже недефективен и имеет
полуиндекс $0$); в точности $\delta(k/2)$ из этих $|\frac{k}{2}|-1$
классов является гео\-мет\-ри\-че\-ски специальным (и даже состоящим из
строго геометри\-чески самосокращающих точек), {\rm см.\ (\ref
{eq:delta}) и определение~\ref {def:2.1}~(b,~c)}.

При четном $k_1k_2$ выполнено свойство Векена, кривые
$\gamma_1,\gamma_2$ гомотопны непересекающимся кривым, а все классы
Нильсена недефективны (и несущественны). При нечетном $k_1k_2$ кривые
$\gamma_1,\gamma_2$ имеют в точности $\gcd(k_1,k_2)$ классов Нильсена
точек пересечения; каждый из этих классов существен и содержит не
менее $\frac{\min\{|k_1|,|k_2|\}}{\gcd(k_1,k_2)}$ точек пересечения,
хотя дефективен.
 \end{unLem}

\begin{unCor}\label{cor:geomMoebius}
В условиях леммы~$\ref {lem:Moebius}$ для замкнутой кривой
$\gamma\colon\, S^1 \to M$ на листе Мебиуса отличные от нуля
полуиндексы классов геометрических точек самопересечения равны $1$ и
   $$
MI_\geom[\gamma] = [\frac{|k|-1}{2}], \qquad RI_\geom[\gamma] =
[\frac{|k|}{2}]+1,
 $$
 $$
 NI_\geom^*[\gamma] =  NI_\geom[\gamma] = \left\{
\begin{array}{ll} \delta(\frac{k}{2}), & \mbox{если $k$ четно}, \\
                        \frac{|k|-1}{2}, & \mbox{если $k$ нечетно}.
\end{array} \right.
   $$

Если $k$ нечетно, то выполнено свойство Векена, все классы Нильсена
геометрических точек самопересечения строго дефективны, а все классы кроме
тривиального существенны.

Если $k$ четно, то кривая $\gamma$ имеет не менее $[|k|/4]$ классов
геометрических точек самопересечения, отличных от три\-ви\-аль\-ного;
в точности $\delta(k/2)$ из этих классов является (строго) геомет\-ри\-чески
специальным (он существен и строго дефективен), а каждый из остальных
$[(|k|-2)/4]$ классов содержит не менее двух (строго эквивалентных)
геометрических точек, хотя несу\-щест\-вен (и даже недефективен и
имеет полуиндекс $0$). \qed
\end{unCor}

\PROOF леммы \ref {lem:Moebius}. Пусть $k$ четно. С
исполь\-зо\-ва\-нием поднятия $\tilde\gamma$ кривой $\gamma$ на
двулистное накрытие $C\to M$ листа Мебиуса $M$ цилиндром $C$, из
леммы~\ref {lem:Cylinder} получаем неравенство $MI[\gamma]\ge
MI[\tilde\gamma]=|k|-2$. Оно является равенством, что доказывается
предъявлением подходящей кривой $\gamma'\simeq\gamma$, см.\ рис.~3.
Так как $\gamma$ сохраняет ориентацию, все классы Нильсена
недефективны. Поскольку $MI[\gamma]=MI[\tilde\gamma]$, из свойств
накрытий и из леммы~\ref {lem:Cylinder} следуют равенства
$NI[\gamma]=NI[\tilde\gamma]=0$; равенство нулю полуиндексов всех
классов Нильсена. Из леммы~\ref {lem:Cylinder} получаем также
существование не менее $|k|/2-1$ нетрививальных классов, каждый из
которых недефективен и содержит не менее двух строго эквивалентных
точек самопересечения; среди них существование в точности $\delta(k)$ геометрически специального класса; существенность этого класса и то, что его точки являются строго геомет\-ри\-чески самосокращающими.

\begin{figure}
\unitlength = 3.8mm
\begin{center}
\begin{picture}(28.75,9.3)(0,-0.5)
\linethickness{2pt} \Thicklines
\path(0,-0.3)(12.5,-0.3)(12.5,7.3)(0,7.3)(0,-0.3)
\path(12.75,0.2)(12.5,-0.3)
\path(12.25,0.2)(12.5,-0.3)
\path(-0.25,6.80)(0,7.3)
\path( 0.25,6.80)(0,7.3)
\put(-0.7,0.8){\small 6}
\put(-0.7,1.8){\small 4}
\put(-0.7,2.8){\small 2}
\put(-0.7,3.8){\small 1}
\put(-0.7,4.8){\small 3}
\put(-0.7,5.8){\small 5}

\put(12.8,0.8){5}
\put(12.8,1.8){3}
\put(12.8,2.8){1}
\put(12.8,3.8){2}
\put(12.8,4.8){4}
\put(12.8,5.8){6}

\thinlines
\path(0,1)(12.5,3)
\path(0,2)(12.5,1)
\path(0,3)(12.5,2)
\qbezier[50](0,3.5)(6.25,3.5)(12.5,3.5)
\put(3.6,3.05){\tiny \sl центральная окружность $\gamma_0$}
\path(0,4)(12.5,4)
\put(10,3.81){$>$}
\path(0,5)(12.5,5)
\path(0,6)(12.5,6)


\Thicklines
\path(16,-0.3)(28.5,-0.3)(28.5,8.3)(16,8.3)(16,-0.3)
\path(28.75,0.2)(28.5,-0.3)
\path(28.25,0.2)(28.5,-0.3)
\path(15.75,7.80)(16,8.3)
\path(16.25,7.80)(16,8.3)
\put(15.3,0.8){\small 6}
\put(15.3,1.8){\small 4}
\put(15.3,2.8){\small 2}
\put(15.3,3.8){\small 1}
\put(15.3,4.8){\small 3}
\put(15.3,5.8){\small 5}
\put(15.3,6.8){\small 7}

\put(28.8,0.8){\small 7}
\put(28.8,1.8){\small 5}
\put(28.8,2.8){\small 3}
\put(28.8,3.8){\small 1}
\put(28.8,4.8){\small 2}
\put(28.8,5.8){\small 4}
\put(28.8,6.8){\small 6}

\thinlines \path(16,1)(28.5,2) \path(16,2)(28.5,3)
\path(16,3)(28.5,4) \put(17,2.92){$<$} 
\path(16,4)(28.5,1) \qbezier[50](16,4)(22.25,4)(28.5,4)
\put(18.5,4.1){\tiny \sl центральная окружность $\gamma_0$}
\path(16,5)(28.5,5) \path(16,6)(28.5,6) \path(16,7)(28.5,7)

\end{picture}
\end{center}

\hskip 1truecm
\begin{tabular}{c} \small
Рис.~3. Специальная кривая \\ $\gamma\simeq\gamma_0^6$
\end{tabular}
\hskip 1truecm
\begin{tabular}{c} \small
Рис.~4. Неспециальная кривая \\ $\gamma\simeq\gamma_0^{-7}$
\end{tabular}
\end{figure}

Пусть $k$ нечетно. Покажем, что $MI[\gamma]=NI[\gamma]=|k|-1$.
Обозначим через $\gamma_0$ центральную окружность листа Мебиуса.
Разъединим концы кривой $\gamma$ и продеформируем ее в простую кривую
$\gamma'$, идущую ``параллельно'' кривой $\gamma_0$ ``по спирали'',
как показано на рис.~4. При такой ``намотке'' кривая $\gamma'$
обходит $|k|>1$ раз вокруг $\gamma_0$ и, следовательно, остается
незамкнутой. ``Замкнем'' кривую $\gamma'$ коротким отрезком $\sigma$,
пересекающим $\gamma_0$ в единственной точке $\gamma_0(v)$; получим
замкнутую кривую $\gamma''=\gamma'\cdot\sigma$, гомотопную $\gamma$,
см.\ рис.~4. По построению (незамкнутые) кривые $\gamma',\sigma$ не
имеют точек самопересечения и пересекаются в $\frac{|k|-1}{2}$
точках, не считая концов. Так как все эти точки определяют точки
самопересечения кривой $\gamma''$, то кривая $\gamma''$ имеет $|k|-1$
точек самопересечения. Покажем, что эти точки попарно неэквивалентны
(а потому не являются геометрически специальными) и нетривиальны.
Сопоставим этим точкам целые числа $g_{v_1,v_2}$ аналогично шагу~2
доказательства леммы~\ref {lem:Cylinder}. Тогда попарная
неэквивалентность следует из того, что $|k|-1$ чисел
$\pm2,\pm4,\dots,\pm(|k|-1)$ попарно различны в факторгруппе
$\ZZ/|k|\ZZ$, а нетривиальность --- из того, что эти числа отличны от
нуля в $\ZZ/|k|\ZZ$. Каждая из этих точек является строго самосокращающей (а
ее класс строго дефективным) ввиду ее специальности и того, что $\gamma$
меняет ориентацию.

Пусть $k_1k_2$ четно. Тогда одно из чисел $k_1,k_2$ четно, скажем,
$k_1$. Значит, кривая $\gamma_1$ гомотопна некоторой степени границы
листа Мебиуса. Продеформируем кривую $\gamma_1$ к границе, а кривую
$\gamma_2$ --- к центральной окружности листа Мебиуса. Получим пару
непе\-ре\-се\-ка\-ю\-щихся кривых, поэтому
$MI[\gamma_1,\gamma_2]=NI[\gamma_1,\gamma_2]=0$. Недефективность
классов Нильсена следует из того, что одна из кривых
$\gamma_1,\gamma_2$ сохраняет ориентацию.

Пусть $k_1k_2$ нечетно. Разобьем доказательство на несколько шагов,
аналогичных шагам доказательства леммы~\ref {lem:Cylinder}.

{\it Шаг 1.} Докажем равенство
$MI[\gamma_1,\gamma_2]=\min\{|k_1|,|k_2|\}$. Рассмот\-рим пару кривых
$(\gamma_1',\gamma_2')$, $\gamma_1'\simeq\gamma_1$,
$\gamma_2'\simeq\gamma_2$, имеющих минимальное число точек
пересечения, т.е.\ $MI[\gamma_1,\gamma_2]$ точек пересечения. Как в
доказательстве леммы~\ref {lem:Cylinder} (шаги 1 и 3), каждую пару
кривых $\gamma_1',\gamma_2'$ на листе Мебиуса, имеющую конечное число
точек пересечения, можно продеформировать --- с сохранением точек
пересечения и их разбиения на классы Нильсена --- так, чтобы одна из
кривых $\gamma_1',\gamma_2'$ имела общую точку с $\partial M$, а
другая кривая лежала внутри $M$. Пусть для определенности
$\gamma_2'(0)\in\partial M$, $\gamma_1'(t)\in M\setminus\partial M$.
Рассмотрим ориентируемое двулистное накрытие $C=\tilde M\to M$ листа
Мебиуса $M$ цилиндром $C$, и пусть $\tilde\gamma_1,\delta_2$ --
поднятия кривых $(\gamma_1')^2,\gamma_2'$ на это накрытие. Тогда
первая кривая $\tilde\gamma_1$ замкнута, лежит во внутренности
$C=\tilde M$ и гомотопна $k_1$-ой степени центральной окружности
цилиндра $C$. Вторая кривая $\delta_2$ не является замкнутой и ее
концы $\delta_2(0),\delta_2(1)$ принадлежат разным компонентам края
цилиндра $C=\tilde M$.

Рассмотрим задачу о минимальном числе $MI[\tilde\gamma_1,\delta_2]$
точек пере\-се\-че\-ния замкнутой кривой $\tilde\gamma_1'$ и дуги
$\delta_2'$ на цилиндре в классах гомотопии замкнутой кривой
$\tilde\gamma_1$ и дуги $\delta_2$ с фиксированными концами (в
процессе гомотопии концы дуги остаются неподвижными). Как и в задаче
о пересечении двух дуг $\delta_1,\delta_2$ (см.\ шаг~1 доказательства
леммы~\ref {lem:Cylinder}), в этой задаче справедливо свойство
Векена, и минимальное число точек пересечения в данных классах
гомотопии имеют кривые, являющиеся кратчайшими относительно некоторой
плоской рима\-но\-вой метрики на цилиндре. Отсюда легко находим
 $$
MI[\tilde\gamma_1,\delta_2]=NI[\tilde\gamma_1,\delta_2]=|k_1|.
 $$
Следовательно, кривые $\tilde\gamma_1,\delta_2$ имеют не менее
$|k_1|$ точек пересечения.

Так как разным точкам пересечения кривых $\tilde\gamma_1,\delta_2$
отвечают разные точки пересечения кривых $\gamma_1',\gamma_2'$,
получаем неравенство $MI[\gamma_1,\gamma_2]\ge
MI[\tilde\gamma_1,\delta_2]=|k_1|\ge\min\{|k_1|,|k_2|\}$. Оно
является равенством $MI[\gamma_1,\gamma_2]=\min\{|k_1|,|k_2|\}$, что
доказывается путем предъявления подходящей пары кривых
$\gamma_1',\gamma_2'$, имеющих ровно $|k_i|=\min\{|k_1|,|k_2|\}$
точек пересечения (например, $\gamma_i'=\gamma_0^{k_i}$, а
$\gamma_{3-i}'$ -- любая кривая, гомотопная $\gamma_{3-i}$ и
пересекающая центральную окружность $\gamma_0$ листа Мебиуса в ровно
одной точке).

{\it Шаг 2.} Для вычисления числа Нильсена $NI[\gamma_1,\gamma_2]$
при нечетном $k_1k_2$ сопоставим каждой точке пересечения $(v_1,v_2)$
кривых $\gamma_1,\gamma_2$ путь
$\gamma_{v_1,v_2}=\gamma_1|_{[0,v_1]}\cdot\gamma_2|_{[0,v_2]}^{-1}$,
начинающийся в точке $\gamma_1(0)$ и заканчивающийся в точке
$\gamma_2(0)$. Гомотопический класс (относительно концов) этого пути
обозначим через
$$
g_{v_1,v_2}=[\gamma_{v_1,v_2}]\in\pi_1(C,\gamma_1(0),\gamma_2(0))\cong\ZZ
$$
и будем рассматривать как целое число. Из первой формулы леммы~\ref
{lem:topol}~(a) следует, что точки пересечения $(v_1,v_2)$ и
$(v_1',v_2')$ эквивалентны в том и только том случае, когда число
$g_{v_1,v_2}-g_{v_1',v_2'}$ является целочисленной линейной
комбинацией чисел $k_1$ и $k_2$. Пусть для простоты $\gamma_1(0)=\gamma_2(0)$. Тогда классы Нильсена пары кривых $\gamma_1,\gamma_2$ находятся во взаимно-однозначном
соответствии с образами чисел $g_{v_1,v_2}$ при проекции группы $\ZZ$
в группу $\ZZ/r\ZZ$, где $r=\gcd(k_1,k_2)$.

Обозначим через $\tilde\gamma_0$ центральную окружность цилиндра.
Рассмот\-рим кратчайшую замкнутую кривую
$\tilde\gamma_1'=\tilde\gamma_0^{k_1}$ и кратчайшую дугу $\delta_2'$
в классах гомотопии $\tilde\gamma_1'\simeq\tilde\gamma_1$,
$\delta_2'\simeq_\partial\delta_2$ относительно какой-нибудь плоской
римановой метрики на цилиндре, где в процессе гомотопии пути
$\delta_2$ его концы $\delta_2(0),\delta_2(1)\in\partial C$
предполагаются неподвижными, см.\ шаг~1. Тогда кривые
$\tilde\gamma_1',\delta_2'$ имеют $|k_1|$ точек пересечения;
отвечающие им $|k_1|$ чисел $g_{v_1,v_2}$ образуют множество
$\{2,4,\dots,2|k_1|\}$. Это множество распадается на $r$ наборов из
$\frac{|k_1|}{r}$ чисел в каждом, где числа каждого набора равны
между собой в группе $\ZZ/r\ZZ$, а элементы разных наборов различны в
$\ZZ/r\ZZ$. Следовательно, для кривых $\gamma_1',\gamma_2'$ число
классов Нильсена точек пересечения равно
$NI[\gamma_1,\gamma_2]=r=\gcd(k_1,k_2)$, и каждый класс Нильсена
состоит из нечетного числа $\frac{|k_1|}{r}$ точек пересечения.
Поэтому все классы Нильсена существенны. Классы Нильсена дефективны,
так как они специальны и обе кривые $\gamma_1,\gamma_2$ меняют
ориентацию.

{\it Шаг 3.} Формулы для $RI[\gamma]$ и $RI[\gamma_1,\gamma_2]$, а
также то, что для любой пары кривых, гомотопных паре
$\gamma_1,\gamma_2$, каждый класс Нильсена точек пересечения содержит
не менее $\frac{\min\{|k_1|,|k_2|\}}{\gcd(k_1,k_2)}$ точек,
доказывается как в доказательстве леммы~\ref {lem:Cylinder} (шаги~3
и~4). \qed \ENDPROOF

\section{СПЕЦИАЛЬНЫЕ КРИВЫЕ И СПЕЦИАЛЬНЫЕ ПАРЫ КРИВЫХ}\label{sec:4}

Вводимое ниже определение мотивировано тем, что в задаче о
(само)пересечении замкнутых кривых на поверхности свойство Векена
выполнено для неспециальных кривых и неспециальных пар привых, но оно
нарушается для специальных кривых и специальных пар кривых, причем за
счет специальных классов Нильсена (см.\ теоремы~\ref {pro:2.4} и \ref
{pro:2.2} ниже).

\begin{unDef}\label{def:para}
(a)\quad Замкнутую кривую $\gamma\colon\,S^1\to S$ на поверхности $S$
назовем {\it специальной}, если $\gamma$ сохраняет ориентацию,
нестягиваема и гомотопна собственной степени некоторой замкнутой
кривой $\gamma_0$ на $S$, т.е.\ $\gamma\simeq\gamma_0^k$ для
некоторого числа $k\in\ZZ$, $|k|>1$. Здесь и далее
$\gamma_0^k(v)=\gamma_0(kv)$, $v \in S^1$.

\smallskip
(b)\quad Пару замкнутых кривых $\gamma_1,\gamma_2\colon\,S^1\to S$ на
поверхности $S$ назовем {\it специальной}, если обе кривые
$\gamma_1,\gamma_2$ меняют ориентацию и гомотопны степеням одной и
той же замкнутой кривой $\gamma$ на $S$, т.е.\
$\gamma_1\simeq\gamma^k$, $\gamma_2\simeq\gamma^\ell$ для некоторых
чисел $k,\ell\in\ZZ$.
\end{unDef}

\begin{unRem}\label{rem:osob}
(a) Связная поверхность $S$ содержит специ\-аль\-ную замкнутую кривую
тогда и только тогда, когда ее фунда\-мен\-таль\-ная группа
бесконечна (т.е.\ когда $\chi(S)\le0$, т.е.\ $S$ не гомеоморфна
плоскости, сфере и проективной плоскости). Если $\chi(S)>0$, то любая
замкнутая кривая на $S$ неспециальна и гомотопна собственной степени
замкнутой кривой.

(b) Если либо $S$ связна и $\chi(S)\ge0$, либо кривая $\gamma$
стягиваема (соответственно одна из кривых $\gamma_1,\gamma_2$
стягиваема), то все точки самопересечения и геометрические точки
самопересечения (соответственно все точки пересечения) являются
специальными, см.\ определение~\ref {def:2.1}~(c).

Для остальных связных поверхностей $S$, $\chi(S)<0$, любая кривая
$\gamma\simeq\gamma_0^k$ имеет не более $|k|$ специальных (а потому
не более $|k|-1$ нетривиальных специальных) классов Нильсена точек
самопересечения, а любая пара кривых
$\gamma_1\simeq\gamma_0^{k_1},\gamma_2\simeq\gamma_0^{k_2}$ имеет не
более $\gcd(k_1,k_2)$ специальных классов Нильсена точек пересечения
(здесь $k,k_1,k_2\in\ZZ$, $k\ne0$, $(k_1,k_2)\ne(0,0)$ и $\gamma_0$
--- замкнутая кривая, не гомотопная собственной степени никакой замкнутой
кривой на $S$). При этом специальные точки (само)пересечения --- в
точности те, которые поднимаются на накрытие поверхности $S$
(гомеоморфное цилиндру или листу Мебиуса), отвечающее циклической
подгруппе $Z_{\gamma_0}\subset\pi_1(S,\gamma_0(0))$, см.\
доказательства теорем~\ref {pro:2.4} и~\ref {pro:2.2}.

(c) Можно показать, что если $S$ --- не бутылка Клейна, то для
специальной точки самопересечения $(v_1,v_2)$ кривые
$\gamma(v_1+\overline{t})$ и $\gamma(v_2+\overline{t})$, $0\le t\le
1$, гомотопны относительно концов; т.е.\ можно считать, что $p=q=1$
(для бутылки Клейна $p=1,q=\pm1$). Если $S$ --- не бутылка Клейна и у
пары замкнутых кривых $\gamma_1,\gamma_2$ на $S$ имеется специальная
точка пересечения, то $\gamma_1\simeq\gamma_0^{k_1}$,
$\gamma_2\simeq\gamma_0^{k_2}$ для некоторых замкнутой кривой
$\gamma_0$ и чисел $k_1,k_2\in\ZZ$.
 \end{unRem}

\section{ЗАДАЧА О САМОПЕРЕСЕЧЕНИИ КРИВЫХ НА ПОВЕРХНОСТЯХ}\label{sec:5}

\begin{unThe}\label{pro:2.4}
Пусть $\gamma\colon\, S^1 \to S$ -- замкнутая кривая на поверхности
$S$, не обязательно являющейся компактной. Тогда отличные от нуля
полуиндексы классов Нильсена точек самопересечения кривой $\gamma$
равны $1$. Более того:

\smallskip
{\rm (a)}\quad Неспециальная кривая $\gamma$ обладает свойством
Векена для задачи о самопересечении:
   $$
MI[\gamma] = NI^*[\gamma] = NI[\gamma]\,.
   $$
Если при этом $\gamma\simeq\gamma_0^k$, где $k\in\ZZ\setminus\{0\}$,
и кривая $\gamma_0$ негомотопна собственной степени никакой замнутой
кривой на $S$, то
   $$
MI[\gamma] = NI^*[\gamma]=NI[\gamma]=k^2\cdot NI[\gamma_0] + |k|-1
   $$
и ровно $|k|-1$ существенных классов Нильсена точек
само\-пе\-ре\-се\-чения кривой $\gamma$ являются специальными, причем
эти классы отличны от тривиального,
строго дефективны и не являются геометрически специальными {\rm (см.\
определение~\ref {def:2.1}(b,~c))}.

\smallskip
{\rm (b)}\quad Специальная кривая $\gamma$ гомотопна степени
$\gamma_0^k$, $k\in\ZZ$, $|k|>1$, где кривая $\gamma_0$ сохраняет
ориентацию или $k$ четно, причем $\gamma_0$ негомотопна собственной
степени никакой замнутой кривой на $S$. Положим $k'=k$, если
$\gamma_0$ сохраняет ориентацию, и $k'=\frac{k}{2}$, если $\gamma_0$
меняет ориентацию. Тогда все классы точек самопересечения
недефективны, все специальные классы несущественны и
   $$
 MI[\gamma] = k^2 \cdot NI[\gamma_0] + 2(|k'|-1) ,\qquad
 NI^*[\gamma] =  NI[\gamma] = k^2 \cdot NI[\gamma_0] ,
   $$
причем $\gamma$ имеет не менее $|k'|-1$ (строго) специальных классов точек
само\-пересечения, каждый из которых не совпадает с триви\-аль\-ным и
содержит не менее двух (строго эквивалентных) точек, хотя несуществен
(и даже недефективен и имеет полу\-индекс $0$); в точности
$\delta(k')$ из этих $|k'|-1$ классов явля\-ется геометрически
специальным (и даже состоящим из строго геомет\-ри\-чески
самосокращающих точек), {\rm см.~(\ref {eq:delta}) и опре\-де\-ле\-ние~\ref
{def:2.1}~(b,~c)}.
 \end{unThe}

Согласно теореме~\ref {pro:2.4}, в задаче о самопересечении замкнутой
кривой $\gamma$ свойство Векена выполняется в точности для
неспециаль\-ных кривых и для специальных кривых вида $\gamma_0^2$,
где кривая $\gamma_0$ меняет ориентацию и негомотопна собственной
степени никакой замкнутой кривой.

\begin{unCor}\label{cor:geomselfint}
В условиях теоремы~$\ref {pro:2.4}$ отличные от нуля полуиндексы
классов Нильсена геометрических точек самопересечения равны $1$;

\smallskip
{\rm (a)}\quad для неспециальной кривой $\gamma\simeq\gamma_0^k$
выполнено свойство Векена
   $$
 MI_\geom[\gamma] =
 NI_\geom^*[\gamma] = NI_\geom[\gamma] = k^2\cdot NI_\geom[\gamma_0]+\frac{|k|-1}{2} 
   $$
и ровно $\frac{|k|-1}{2}$ существенных классов специальны, причем эти
классы строго специальны, строго дефективны и отличны от тривиального;

\smallskip
{\rm (b)}\quad для специальной кривой $\gamma\simeq\gamma_0^k$
выполнено
   $$
\begin{array}{rcl}
MI_\geom[\gamma] &=& k^2 \cdot NI_\geom[\gamma_0] + |k'|-1 ,\\
NI_\geom^*[\gamma] = NI_\geom[\gamma] &=& k^2 \cdot
NI_\geom[\gamma_0] + \delta(k') 
\end{array}
   $$
и ровно $k^2 \cdot NI_\geom[\gamma_0]$ существенных классов неспециальны,
где $k'=k$, если $\gamma_0$ сохраняет ориентацию, и $k'=\frac{k}{2}$,
если $\gamma_0$ меняет ориентацию, причем кривая $\gamma$ имеет не
менее $[|k'|/2]$ (строго) специальных классов геометрических точек самопересечения,
отличных от три\-ви\-аль\-ного; в точности $\delta(k')$ из этих
классов является (строго) геомет\-ри\-чески специальным (он существен и строго
дефективен), а каждый из остальных $[(|k'|-1)/2]$ классов содержит не
менее двух (строго эквивалентных) гео\-мет\-ри\-чес\-ких точек, хотя
несу\-щест\-вен (и даже недефективен и имеет полуиндекс $0$).
 \qed
\end{unCor}

Согласно следствию \ref {cor:geomselfint}, свойство Векена для
геометрических точек самопересечения выполнено в точности для
неспециальных кривых, а также для специальных кривых вида
$\gamma_0^4$ ($\gamma_0$ меняет ориентацию) и $\gamma_0^2$, где
$\gamma_0$ не гомотопна степени никакой замкнутой кривой на $S$.

\PROOF теоремы \ref {pro:2.4}. Если $S=\RR P^2$ или кривая $\gamma$
стягиваема, то она неспециальна, и требуемое свойство
$MI[\gamma]=NI[\gamma]$ для этой кривой очевидно. Поэтому будем
считать, что $S$ связна, $|\pi_1(S)|=\infty$ и кривая $\gamma$
нестягиваема. Так как образ любой гомотопии кривой в $S$ является
компактным подмножеством в $S$, можно считать, что $S$ является
компактной поверхностью (с краем или без) и $\chi(S)\le0$. Известно,
что в этом случае на $S$ существует риманова метрика, у которой
кривизна всюду неположительна и каждая компонента края является
геодезической.

В силу допущенных предположений, существует замкнутая кривая
$\gamma_0$, не гомотопная собственной степени никакой замкнутой
кривой на $S$, такая, что $\gamma\simeq\gamma_0^k$, $k\in\ZZ$, $k>0$.
Пусть $\gamma_0'$ -- кратчайшая кривая, гомотопная кривой $\gamma_0$
в $S$; тогда $\gamma_0'$ является замкнутой геодезической (и
наоборот: любая замкнутая геодезическая является кратчайшей замкнутой
кривой в своем классе гомотопии). При этом в любой точке
самопересечения геодезической $\gamma_0'$ ее ветви пересекаются
трансверсально (иначе кривая $\gamma_0$ являлась бы собст\-вен\-ной
степенью другой геодезической). Далее будем обозначать кривую
$\gamma_0'$ через $\gamma_0$, и будем считать, что
$\gamma=\gamma_0^k$.

{\it Случай 1.} Предположим, что $k=1$, т.е.\ $\gamma=\gamma_0$. Если
бы у $\gamma$ существовала специальная точка самопересечения или две
эквивалентные точки самопересечения (например, эквивалентные точки
вида $(v_1,v_2)$ и $(v_2,v_1)$), то существовали бы два поднятия
$\tilde\gamma_1,\tilde\gamma_2$ кривой $\gamma$ в универсальную
накрывающую поверхность $\tilde S$, имеющие две точки
трансверсального пересечения, либо одно из поднятий
$\tilde\gamma_1,\tilde\gamma_2$ являлось бы замкнутой кривой. Это
противоречит тому факту, известному еще Пуанкаре~\cite {Po}, что в
полном односвязном римановом многообразии неположительной секционной
кривизны все геодезические незамкнуты и любые две из них либо
совпадают, либо имеют не более одной общей точки.

Следовательно, точки самопересечения кривой $\gamma$ неспециальны,
попарно неэквивалентны (в том числе точки вида $(v_1,v_2)$ и
$(v_2,v_1)$) и имеют ненулевые индексы $\pm1$, а значит, все классы
Нильсена существенны, нетривиальны и имеют полуиндекс $1$. В
частности, геодезическая $\gamma$ имеет ровно $NI[\gamma]$ точек
самопересечения (точки вида $(v_1,v_2)$ и $(v_2,v_1)$ считаются
разными), откуда $MI[\gamma]\le NI[\gamma]$. Следовательно,
$MI[\gamma]=NI[\gamma]$ и среди существенных классов нет специальных
(а потому нет тривиального).

{\it Случай 2.} Предположим, что $k>1$. Пусть число точек
самопересечения кривой $\gamma_0$ равно $2n$ (точки вида $(v_1,v_2)$
и $(v_2,v_1)$ считаются различными). Из доказательства случая~1
следует, что
 $$
MI[\gamma_0]=NI[\gamma_0]=2n\, .
 $$

Как в доказательствах лемм~\ref {lem:Cylinder} и~\ref {lem:Moebius},
разъединим концы кривой $\gamma$ и продеформируем ее в простую кривую
$\gamma'$, которая идет ``параллельно'' кривой $\gamma_0$ ``по
спирали'' (см.\ рис.~3,~4). При такой ``намотке'' кривая $\gamma'$
обходит $k>1$ раз вокруг $\gamma_0$ и, следовательно, остается
незамкнутой. ``Замкнем'' кривую $\gamma'$ коротким отрезком $\sigma$;
получим замкнутую кривую $\gamma''=\gamma'\cdot\sigma$, гомотопную
$\gamma$.

По построению (незамкнутая) кривая $\gamma'$ имеет $k^2\cdot
NI[\gamma_0]$ точек самопересечения. Если бы одна из этих точек
являлась специальной точкой самопересечения кривой $\gamma''$ или две
из этих точек были эквивалентны, то геодезическая $\gamma_0$ имела бы
либо специальную точку самопересечения, либо тривиальную точку
самопересечения, либо две эквивалентные точки самопересечения, что
противоречит доказанному в случае~1. По построению (незамкнутая)
кривая $\sigma$ не имеет точек самопересечения и пересекает
(незамкнутую) кривую $\gamma'$ в
 $$
N:=\left\{ \begin{array}{ll}
 \frac{k-1}{2}, & \mbox{если $\gamma$ -- неспециальная}, \\
          k'-1, & \mbox{если $\gamma$ -- специальная}, \\
 \end{array} \right.
 $$
точках, причем все эти точки определяют специальные точки
самопе\-ресечения кривой $\gamma''=\gamma'\cdot\sigma$.
Следовательно, кривая $\gamma''$ имеет $k^2\cdot NI[\gamma_0]+2N$
точек самопересечения, из которых $k^2\cdot NI[\gamma_0]$ точек
неспециальны и попарно неэквивалентны, а остальные $2N$ точек
специальны.

В случае (a) осталось показать, что специальные точки кривой
$\gamma''$ тоже попарно неэквивалентны, а также нетривиальны и строго
дефективны. Отсюда будет, в частности, следовать, что среди
сущест\-вен\-ных классов Нильсена имеется ровно $k^2\cdot
NI[\gamma_0]$ неспециальных и $2N=k-1$ специальных классов, причем
нет тривиального класса и что $MI[\gamma]=NI[\gamma]=k^2\cdot
NI[\gamma_0]+k-1$. В случае (b) осталось показать, что любая кривая,
гомотопная $\gamma$, имеет не менее $k'-1$ нетривиальных специальных
классов Нильсена точек самопересечения, каждый из которых содержит не
менее двух (строго эквивалентных) точек, недефективен и имеет
полу\-индекс $0$, причем в точности $\delta(k')$ из этих $|k'|-1$
классов состоит из строго геомет\-ри\-чески самосокращающих точек.

Пусть $*=\gamma_0(0)$ -- базисная точка поверхности $S$,
$[\gamma_0]\in\pi_1(S,*)$ -- гомотопический класс кривой $\gamma_0$ в
$S$, и пусть $p\colon\,\tilde S\to S$ -- накрытие поверхности $S$ с
базисной точкой $\tilde*\in p^{-1}(*)$, отвечающее циклической
подгруппе $p_\#\pi_1(\tilde S,\tilde*)=\langle[\gamma_0]\rangle$
фундаментальной группы $\pi_1(S,*)$ с образующей $[\gamma_0]$.
Обозначим через $P$ цилиндр, если $\gamma_0$ сохраняет ориентацию,
или лист Мебиуса, если $\gamma_0$ меняет ориентацию. Поскольку
$\chi(S)\le0$, в $\pi_1(S,*)$ нет элементов конечного порядка;
значит, $\pi_1(\tilde S,\tilde*)\approx\ZZ$ и $\tilde S$ гомеоморфно
поверхности $P\setminus B$ без некоторого подмножества края
$B\subset\partial P$, причем поднятие $\tilde\gamma_0$ кривой
$\gamma_0$ в $\tilde S$ является центральной окружностью $P$. Так как
кривая $\gamma_0$ допускает поднятие $\tilde\gamma_0$, то гомотопия
$\gamma_0^k\simeq\gamma''$ допускает поднятие
$(\tilde\gamma_0)^k\simeq\tilde\gamma''$. Ясно, что если точка
самопересечения кривой $\gamma''$ ``поднимается'', то она является
специальной и все эквивалентные ей точки самопересечения тоже
``поднимаются'' и определяют специальные точки самопересечения кривой
$\tilde\gamma''$. Как следует из свойств накрытий, тривиальность (или
эквивалентность) этих точек по отношению к ``поднятой'' кривой
$\tilde\gamma''$ равносильна их тривиальности (соответственно
эквивалентности) по отношению к самой кривой $\gamma''$.

(a) Поэтому достаточно показать, что точки самопересечения
``поднятой'' кривой $\tilde\gamma''$ на листе Мебиуса попарно
неэквивалентны, а также нетривиальны и строго дефективны --- это
следует из леммы~\ref {lem:Moebius}. Поэтому
 $$
MI[\gamma]=NI[\gamma]=k^2\cdot NI[\gamma_0]+k-1.
 $$

(b) Достаточно показать, что любая кривая на $P$, гомотопная
$\tilde\gamma''$,
имеет не менее $k'-1$ нетривиальных классов Нильсена точек
самопересечения, каждый из которых содержит не менее двух (строго
эквивалентных) точек, недефективен и имеет полу\-индекс $0$, причем в
точности $\delta(k')$ из этих $|k'|-1$ классов состоит из строго
геомет\-ри\-чески самосокращающих точек. Это следует из лемм~\ref
{lem:Cylinder} и~\ref {lem:Moebius} для сохраняющей ориентацию кривой
$\gamma$.
 \qed
\ENDPROOF

\section{ЗАДАЧА О ПЕРЕСЕЧЕНИИ ПАР КРИВЫХ НА ПОВЕРХНОСТЯХ}\label{sec:6}

\begin{unThe}\label{pro:2.2}
Пусть $\gamma_1,\gamma_2\colon\, S^1 \to S$ -- две замкнутые кривые
на поверхности $S$, которая не обязательно компактна. Отличные от
нуля полуиндексы классов Нильсена точек пере\-се\-чения равны $1$.

\smallskip
{\rm (a)} Пусть $(\gamma_1,\gamma_2)$ --- неспециальная пара. Тогда
все сущест\-венные классы Нильсена точек пересечения кривых
$\gamma_1,\gamma_2$ явля\-ются неспециальными и выполняется свойство
Векена для задачи о пересечении
   $$
MI[\gamma_1,\gamma_2] = NI^*[\gamma_1,\gamma_2] =
NI[\gamma_1,\gamma_2]\,.
   $$
Если $\gamma_1\simeq\gamma^k$, $\gamma_2\simeq\gamma^\ell$, где
$k,\ell\in\ZZ$ и кривая $\gamma$ негомотопна собственной степени
никакой замкнутой кривой на $S$, то
   $$
MI[\gamma_1,\gamma_2] =
NI^*[\gamma_1,\gamma_2]=NI[\gamma_1,\gamma_2]=|k\cdot\ell|\cdot
NI[\gamma] \, .
   $$

\smallskip
{\rm (b)} Пусть $(\gamma_1,\gamma_2)$ --- специальная пара. Если
$S=\RR P^2$, то $\gamma_1\simeq\gamma_2\not\simeq0$,
$$
MI[\gamma_1,\gamma_2]=NI^*[\gamma_1,\gamma_2]=NI[\gamma_1,\gamma_2]=1
$$
и все точки пересечения образуют существенный специальный (и даже дефек\-тив\-ный)
класс Нильсена. Если $S\ne\RR P^2$, то $\gamma_1\simeq\gamma^k$,
$\gamma_2\simeq\gamma^\ell$, где числа $k,\ell\in\ZZ$ нечетны, а
замкнутая кривая $\gamma$ меняет ориентацию и негомотопна собственной
степени никакой замнутой кривой на $S$; при этом
   $$
\begin{array}{rcl}
MI[\gamma_1,\gamma_2] &=& |k\cdot\ell|\cdot NI[\gamma] + \min\{|k|,|\ell|\}\, ,\\
NI^*[\gamma_1,\gamma_2] = NI[\gamma_1,\gamma_2] &=& |k\cdot\ell|\cdot
NI[\gamma] + \gcd(k,\ell)\,
\end{array}
   $$
и ровно $\gcd(k,\ell)$ существенных классов точек пересечения кривых
$\gamma_1,\gamma_2$ являются специальными; каждый из этих классов
содержит не менее
$\frac{\min\{|k_1|,|k_2|\}}{\gcd(k_1,k_2)}$ точек пересечения, хотя
дефективен.
\end{unThe}

Согласно теореме~\ref {pro:2.2}, в задаче о пересечении замкнутых
кривых свойство Векена нарушается в точности для тех специальных пар
кривых, где никакая из кривых негомотопна степени другой кривой,
причем нарушение свойства Векена происходит за счет специальных
существенных классов Нильсена.

\PROOFp Случай, когда одна из кривых $\gamma_1,\gamma_2$ стягиваема,
очевиден: пара неспециальна и
$MI[\gamma_1,\gamma_2]=NI[\gamma_1,\gamma_2]=0$. Поэтому можно
считать, что $\pi_1(S)\ne0$ и что обе кривые нестягиваемы. Если
$S=\RR P^2$ и $\gamma_1\simeq\gamma_2\not\simeq0$, то пара является
специальной и равенство
$MI[\gamma_1,\gamma_2]=NI[\gamma_1,\gamma_2]=RI[\gamma_1,\gamma_2]=1$
очевидно; кроме того, легко проверяется, что все точки пересечения
дефективны и попарно эквивалентны по Нильсену. Поэтому далее считаем,
что $|\pi_1(S)|=\infty$ и что обе кривые нестягиваемы.

Так как образ любой гомотопии является компактным связным
подмножеством в $S$, можно считать, что $S$ является компактной
связной поверхностью (с краем или без) и, в силу допущенных
предположений, $\chi(S)\le0$. Известно, что в этом случае на $S$
существует риманова метрика, у которой кривизна всюду
неположи\-тельна и каждая компонента края является геодезической.
Пусть $\gamma_1',\gamma_2'$ -- кратчайшие кривые, гомотопные кривым
$\gamma_1,\gamma_2$ в $S$; тогда $\gamma_1',\gamma_2'$ являются
замкнутыми геодезическими (и наоборот: любая замкнутая геодезическая
является кратчайшей замкнутой кривой в своем классе гомотопии). Далее
будем обозначать $\gamma_1',\gamma_2'$ через $\gamma_1,\gamma_2$
соответственно. Рассмотрим два случая:

{\it Случай 1.} Предположим, что не существует замкнутой кривой
$\gamma$ на $S$ и чисел $k,\ell\in\ZZ$, таких, что
$\gamma_1\simeq\gamma^k$, $\gamma_2\simeq\gamma^\ell$. Тогда пара
$(\gamma_1,\gamma_2)$ неспециальна и геодезические
$\gamma_1,\gamma_2$ пересекаются трансверсально в любой точке
пересечения. Если бы существовала специальная точка пересечения или
две эквивалентные точки пересечения, то поднятия
$\tilde\gamma_1,\tilde\gamma_2$ кривых $\gamma_1,\gamma_2$ в
универсальную накрывающую поверхность $\tilde S$ имели бы две точки
трансверсального пересечения. Последнее противоречит тому факту, что
любые две геодезические в односвязном римановом многообразии
неположительной секционной кривизны либо совпадают, либо имеют не
более одной общей точки. Следо\-ва\-тельно,
$MI[\gamma_1,\gamma_2]=NI[\gamma_1,\gamma_2]$, все существенные
классы Нильсена неспециальны и имеют полуиндекс $1$.

{\it Случай 2.} Предположим, что $\gamma_1\simeq\gamma^k$,
$\gamma_2\simeq\gamma^\ell$ для некоторой замкнутой кривой $\gamma$
на $S$ и чисел $k,\ell\in\ZZ$. Можно считать, что $\gamma$ является
геодезической и не гомотопна собственной степени никакой замкнутой
кривой на $S$, и $\gamma_1=\gamma^k$, $\gamma_2=\gamma^\ell$. Пусть
число точек самопересечения геодезической $\gamma$ равно $2n$ (точки
вида $(v_1,v_2)$ и $(v_2,v_1)$ считаются различными). Из
доказательства теоремы~\ref {pro:2.4}~(a) следует, что
 $$
MI[\gamma]=NI[\gamma]=2n
 $$
и все точки самопересечения кривой $\gamma$ неспециальны и попарно
неэквивалентны (в частности, неэквивалентны точки вида $(v_1,v_2)$ и
$(v_2,v_1)$).

(a) Пусть $(\gamma_1,\gamma_2)$ -- неспециальная пара. Тогда либо
кривая $\gamma$ сохраняет ориентацию, либо одно из чисел $k,\ell$
четно, скажем $\ell$. Поэтому можно сдвинуть кривую $\gamma_2$ в
``нормальном'' к $\gamma$ направлении, так, чтобы новая замкнутая
кривая $\gamma_2'$ шла ``парал\-лельно'' $\gamma$. Тогда $\gamma_1$ и
$\gamma_2'$ имеют $|k\ell|\cdot 2n$ точек пересечения. Если бы
$\gamma_1$ и $\gamma_2'$ имели специальную точку пересечения или две
различные эквивалентные точки пересечения, то кривая $\gamma$ имела
бы либо специальную точку самопересечения, либо две эквивалентные
точки самопересечения (например, эквивалентные точки вида $(v_1,v_2)$
и $(v_2,v_1)$), противоречие. Следовательно, все точки пересечения
кривых $\gamma_1,\gamma_2'$ неспециальны, попарно неэквивалентны и
имеют полуиндексы $1$. Значит,
 $$
MI[\gamma_1,\gamma_2]=NI[\gamma_1,\gamma_2]=|k\cdot\ell|\cdot
NI[\gamma]
 $$
и среди существенных классов нет специальных.

(b) Пусть $(\gamma_1,\gamma_2)$ -- специальная пара. Тогда $\gamma$
меняет ориен\-та\-цию, а числа $k,\ell$ нечетны. Пусть для
определенности $|k|\le|\ell|$.

Как в доказательстве леммы~\ref {lem:Moebius} и теоремы~\ref
{pro:2.4}, разъединим концы кривой $\gamma_2$ и сдвинем ее в
``нормальном'' к $\gamma$ направлении так, чтобы новая (незамкнутая)
кривая $\gamma_2'$ шла ``параллельно'' $\gamma$, при этом начальной
точкой на кривой $\gamma$ выбирается точка $\gamma(v)$, являющаяся
образом лишь одного значения $v\in S^1$ при отображении $\gamma$. При
таком ``сдвиге'' кривая $\gamma_2'$ обходит вдоль $\gamma$ нечетное
число раз $|\ell|$ и, следовательно, остается незамкнутой.
``Замкнем'' кривую $\gamma_2'$ коротким отрезком $\sigma$,
пересекающим $\gamma$ в единственной точке $\gamma(v)$; получим
замкнутую кривую $\gamma_2''=\gamma_2'\cdot\sigma$.

По построению кривые $\gamma_1$ и $\gamma_2'$ имеют
$|k\cdot\ell|\cdot 2n$ точек пересечения. Если бы одна из этих точек
являлась специальной точкой пересечения кривых $\gamma_1,\gamma_2''$
или две из этих точек были эквивалентны, то геодезическая $\gamma$
имела бы либо специальную точку самопересечения, либо две
эквивалентные точки самопересечения (например, две эквивалентные
точки вида $(v_1,v_2)$ и $(v_2,v_1)$), про\-ти\-во\-речие. По
построению кривые $\gamma_1$ и $\sigma$ имеют $|k|$ точек
пересечения, причем все они являются специальными (и даже
самосокращающими, см.\ определение~\ref {def:2.1}~(c)) точками
пересечения кривых $\gamma_1,\gamma_2''$. Таким образом, кривые
$\gamma_1$ и $\gamma_2''=\gamma_2'\cdot\sigma$ имеют ровно
$|k\cdot\ell|\cdot 2n+|k|$ точек пересечения, из которых
$|k\cdot\ell|\cdot 2n$ точек неспециальны и попарно неэквивалентны по
Нильсену, а остальные $|k|=\min\{|k|,|\ell|\}$ точек специальны.

Ниже мы покажем, что кривые $\gamma_1,\gamma_2''$ имеют ровно
$\gcd(k,\ell)$ специальных существенных классов Нильсена точек
пересечения, причем любая пара кривых, гомотопная
$(\gamma_1,\gamma_2'')$, имеет не менее $\frac{|k|}{\gcd(k,\ell)}$
точек пересечения в каждом специальном существенном классе Нильсена.
Этим будут доказаны формулы для $MI[\gamma_1,\gamma_2]$ и
$NI[\gamma_1,\gamma_2]$, а также все остальные требуемые утверждения.

Пусть $*=\gamma(0)$ -- базисная точка поверхности $S$,
$[\gamma]\in\pi_1(S,*)$ -- гомотопический класс кривой $\gamma$ в
$S$, и пусть $p\colon\,\tilde S\to S$ -- накрытие поверхности $S$ с
базисной точкой $\tilde*\in p^{-1}(*)$, отвечающее циклической
подгруппе $p_\#\pi_1(\tilde S,\tilde*)=\langle[\gamma]\rangle$
фундаментальной группы $\pi_1(S,*)$ с образующей $[\gamma]$.
Обозначим через $\tilde\gamma$ поднятие кривой $\gamma$ в $\tilde S$.
Заметим, что $\tilde S$ гомеоморфно листу Мебиуса $M\setminus B$ без
некоторого подмножества его края $B\subset\partial M$, причем кривая
$\tilde\gamma$ является центральной окружностью листа Мебиуса $M$.

Рассмотрим поднятия $\tilde\gamma_1,\tilde\gamma_2''$ кривых
$\gamma_1,\gamma_2''$ на накрывающую поверхность $\tilde S=M\setminus
B$, лежащую в листе Мебиуса $M$. Из леммы~\ref {lem:Moebius} следует,
что кривые $\tilde\gamma_1,\tilde\gamma_2''$ имеют ровно
$\gcd(k,\ell)$ существенных (и дефективных) классов Нильсена точек
пересечения. Заметим, что если точка пересечения кривых
$\gamma_1,\gamma_2''$ является (самосокращающей) точкой пересечения
``поднятых'' кривых $\tilde\gamma_1,\tilde\gamma_2''$, то она
является специальной (и самосокращающей) и все эквивалентные ей точки
пересечения тоже ``поднимаются'' и имеют такие же индексы, причем
эквивалентность этих точек по отношению к ``поднятым'' кривым
$\tilde\gamma_1,\tilde\gamma_2''$ равносильна их эквивалентности по
отношению к самим кривым $\gamma_1,\gamma_2''$, откуда полуиндекс
любого класса Нильсена для пары $(\tilde\gamma_1,\tilde\gamma_2'')$
совпадает с его полуиндексом для пары $(\gamma_1,\gamma_2'')$.
Сле\-до\-ва\-тельно, $\gamma_1,\gamma_2''$ имеют не менее
$\gcd(k,\ell)$ специальных существенных (и дефективных) классов
Нильсена точек пересечения. По свойствам накрытий любые кривые
$\gamma_1',\gamma_2'''$, гомотопные кривым $\gamma_1,\gamma_2''$ в
$S$, допускают поднятие в $\tilde S$; при этом ``поднятые'' кривые
$\tilde\gamma_1',\tilde\gamma_2'''$ гомо\-топ\-ны кривым
$\tilde\gamma_1,\tilde\gamma_2''$ в $\tilde S$. По лемме~\ref
{lem:Moebius} любая пара кривых, гомо\-топных
$(\tilde\gamma_1,\tilde\gamma_2'')$ на листе Мебиуса $M$, имеет не
менее $\frac{\min\{|k|,|\ell|\}}{\gcd(k,\ell)}$ точек пересечения в
каждом существенном классе Нильсена. В частности, ``поднятая'' пара
$(\tilde\gamma_1',\tilde\gamma_2''')$ на $\tilde S$ имеет не менее
$\frac{\min\{|k|,|\ell|\}}{\gcd(k,\ell)}$ точек пере\-се\-чения в
каждом существенном классе Нильсена. Следовательно, исходная пара
кривых $(\gamma_1',\gamma_2''')$ на $S$ имеет не менее $\gcd(k,\ell)$
специальных существенных классов Нильсена точек пересечения, каждый
из которых содержит не менее
$\frac{\min\{|k|,|\ell|\}}{\gcd(k,\ell)}$ точек пересечения,
дефективен и имеет полуиндекс $1$.

Как показано выше, среди существенных классов Нильсена имеется ровно
$|k\cdot\ell|\cdot NI[\gamma]$ неспециальных классов.
 \qed
 \ENDPROOF

\REFERENCES{0000}

\bibitem {BKZ} {\it  С.А.\ Богатый, Е.А.\ Кудрявцева, Х.\ Цишанг},
О точках совпадения отображений тора в поверхность, Тр.\ Мат.\ ин-та
им.\ В.\,А.\ Стеклова {\bf 247} (2004), 15--34.

\bibitem{Ni} {\it J.\ Nielsen}, Untersuchungen zur Topologie der geschlossenen
zweiseitigen Fl\"achen I, II, III, Acta Math. {\bf 50} (1927),
189-358; {\bf 53} (1929), 1-76; {\bf 58} (1932), 87-167.

\bibitem{DobJez} {\it R.\ Dobre\'nko, J.\ Jezierski},
The coincidence Nielsen number on nonorientable manifolds, Rocky
Mount.\ J.\ Math. {\bf 23} (1993), 67-87.

\bibitem {Jez} {\it J.\ Jezierski}, The semi-index product formula,
Fund.\ Math. {\bf 140} (1992), 99-120.

\bibitem{We} {\it F. Wecken}, Fixpunktklassen, I, II, III, Math.\ Ann.
{\bf 117} (1941), 659-671; {\bf 118} (1942), 216-234; {\bf 118}
(1942), 544-577.

\bibitem {BGKZ} {\it  С.А.\ Богатый, Д.Л.\ Гонсалвес, Е.А.\ Кудрявцева, Х.\
Цишанг}, Минимальное число прообразов при отображениях между
поверхностями, Мат.\ Заметки {\bf 75}:1 (2004), 13-19.

\bibitem {BGKZ1} {\it S.\ Bogatyi, D.L.\ Gon\c calves, E.\ Kudryavtseva,
H.\ Zieschang}, On the Wecken property for the root problem of
mappings between surfaces, Moscow Math.\ J. {\bf 3}, no.~4 (2003),
1223-1245.

\bibitem {GKZ} {\it D.L.\ Gon\c calves, E.A.\ Kudryavtseva, H.\ Zieschang},
Roots of mappings on nonorientable surfaces and equations in free
groups, Manuscr.\ math. {\bf 107}, no.~3 (2002), 311-341.

\bibitem {Kukh} {\it R.\ Dobre\'nko, Z.\ Kucharski},
On the generalization of the Nielsen number,
Fund.\ Math. {\bf 134}, no.~1 (1990), 1–14.

\bibitem{McCord} {\it C.\ McCord}, A Nielsen theory for intersection numbers,
Fund.\ Math. {\bf 152} (1997), 117-150.

\bibitem {Kervaire} {\it M.A.\ Kervaire}, Geometric and algebraic intersection
numbers, Math.\ Ann. {\bf 103} (1930), 347-358.

\bibitem{Wall} {\it C.T.C.\ Wall}, Surgery on compact manifolds.
London: Academic Press 1970.

\bibitem{Milnor} {\it J.\ Milnor}, Lectures on the $h$-cobordism
theorem, Princeton, New Jersey, Princeton University Press, 1965.

\bibitem {Bo1} {\it B.\ Jiang}, Fixed points and braids,
Invent.\ Math. {\bf 75} (1984), 69-74.

\bibitem {Bo2} {\it B.\ Jiang}, Fixed points and braids II,
Math.\ Ann. {\bf 272} (1985), 249-256.

\bibitem{Ho} {\it H.\ Hopf}, Zur Topologie der Abbildungen von
Mannigfaltigkeiten II, Math.\ Ann. {\bf Bd.~102}  (1930), 562-623.

\bibitem {Bo0} {\it B.\ Jiang}, Fixed points of surface homeomorphisms,
Bull.\ Amer.\ Math.\ Soc. {\bf 5} (1981), 176-178.

\bibitem{Iv} {\it Н.В.\ Иванов}, Числа Нильсена отображений поверхностей,
Зап.\ науч.\ сем.\ ЛОМИ {\bf 122} (1982), 56-65.

\bibitem {Jez1} {\it J.\ Jezierski}, The coincidence Nielsen number for maps
into the real projective space, Fund.\ Math. {\bf 140} (1992),
121-136.

\bibitem {Jez2} {\it J.\ Jezierski}, The least number of coincidence points on
surfaces, J.\ Austr.\ Math.\ Soc.\ A. {\bf 58} (1995), 27-38.

\bibitem {BBPT} {\it R.B.S.\ Brooks, R.F.\ Brown, J.\ Pak, D.H.\ Taylor},
Nielsen numbers of maps of tori, Proc.\ Amer.\ Math.\ Soc. {\bf 52}
(1975), 398-400.

\bibitem{Whitney} {\it H.\ Whitney}, The selfintersection of a smooth
$n$-manifold in $2n$-space, Ann.\ Math. {\bf 45} (1944), 220-246.

\bibitem{Sark} {\it K.S.\ Sarkaria}, Kuratowski complexes,
Topology {\bf 30} (1991), 67-76.

\bibitem {KeMi} {\it M.A.\ Kervaire, J.W.\ Milnor}, On $2$-spheres in
$4$-manifolds, Proc.\ Nat.\ Acad.\ Sci. {\bf 47} (1961), 1651-1657.

\bibitem {Po} {\it H.\ Poincar\'e}, Cinqui\`eme complement \`a l'analysis
situs, Oevres de Henri Poincar\'e, vol.\ VI, 435-498.

\bibitem {TuVi} {\it В.Г.\ Тураев, О.Я.\ Виро}, Пересечения петель в двумерных
многообразиях. II. Свободные петли, Мат.\ Сборник, {\bf 121(163)}
(1983), 350-369; англ.\ пер. {\it V.G.\ Turaev, O.Ya.\ Viro},
Intersections of loops in two-dimensional manifolds. II. Free loops,
Math.\ USSR Sb. {\bf 49} (1984), 357-366.

\bibitem {Ca} {\it P.\ Cahn}, A generalization of the Turaev
cobracket and the minimal self-intersection number of a curve on a
surface. arXiv:1004.0532v3 [math.GT]

\bibitem {CaCh} {\it P.\ Cahn, V.\ Chernov}, Intersection of loops
and the Andersen-Mattes-Reshetikhin algebra. arXiv:1105.4638v1
[math.GT]

\bibitem {Rein} {\it B.\ Reinhardt}, The winding number on two-manifolds,
Ann.\ Inst.\ Fourier Grenoble {\bf 10} (1960), 271-283.

\bibitem {Z} {\it H.\ Zieschang}, Algorithmen f\"ur einfache Kurven auf
Fl\"achen, Math.\ Scand. {\bf 17} (1965), 17-40; II, Math.\ Scand.
{\bf 25} (1969), 49-58.

\bibitem {Ch} {\it D.R.J.\ Chillingworth}, Winding numbers on surfaces. II,
Math.\ Ann. {\bf 199} (1972), 131-153.

\bibitem {BiSe} {\it J.S.\ Birman, C.\ Series}, An algorithm for simple
closed curves on surfaces, J.\ London Math.\ Soc. (2) {\bf 29}
(1984), 331-342.

\bibitem {Lu1} {\it M.\ Cohen, M.\ Lustig}, Paths of geodesics and geometric
intersection numbers. In: Combinatorial group theory and topology
(ed.\ S.M.\ Gersten and J.R.\ Stallings), Ann.\ Math.\ Studies {\bf
111} (1987), 479-500.

\bibitem {Lu2} {\it M.\ Lustig}, Paths of geodesics and geometric
intersection numbers: II. In: Combinatorial group theory and topology
(ed.\ S.M.\ Gersten and J.R.\ Stallings), Ann.\ Math.\ Studies {\bf
111} (1987), 501-543.

\bibitem {Turaev} {\it В.Г.\ Тураев}, Пересечения петель в двумерных
многообразиях, Мат.\ Сборник, {\bf 106(148)} (1978), 566-588; англ.\
пер. Math.\ USSR Sb. {\bf 35} (1979), 229-250.

\bibitem {GKZ2} {\it D.L.\ Gon\c calves, E.A.\ Kudryavtseva, H.\ Zieschang},
Intersection index of curves on surfaces and applications to
quadratic equations in free groups, Atti Sem.\ Mat.\ Fis.\ Univ.\
Modena, Supp.\ {\bf IL} (2001), 339-400.

\bibitem {GKZ3} {\it D.L.\ Gon\c calves, E.A.\ Kudryavtseva, H.\ Zieschang}, An algorithm for
minimal number of intersection points of curves on surfaces, Proc.\
Seminar on Vector and Tensor Analysis {\bf 26} (2005), 139-167.


\bibitem {ScTe} {\it B.\ Schneiderman, P.\ Teichner}, Higher order intersection
numbers of $2$-spheres in $4$-manifolds, Algebraic and Geometric
Topology {\bf 1} (2000), 1-29. http://xxx.lanl.gov/abs/math.GT/0008048

\end{thebibliography}

\end{document}